\documentclass[11pt,a4paper]{amsart}
\usepackage{amsmath,amsthm,amssymb,latexsym}
\usepackage{graphicx}

\newcommand{\bd}{{\mathbb{D}}}
\newcommand{\bi}{{\mathbb{I}}}
\newcommand{\bn}{{\mathbb{N}}}
\newcommand{\br}{{\mathbb{R}}}
\newcommand{\bp}{{\mathbb{P}}}
\newcommand{\bz}{{\mathbb{Z}}}
\newcommand{\bc}{{\mathbb{C}}}
\newcommand{\bo}{{\mathbb{O}}}

\newcommand{\bt}{{\mathbb{T}}}

\newcommand{\ca}{{\mathcal{A}}}

\newcommand{\cd}{{\mathcal{D}}}
\newcommand{\cf}{{\mathcal{F}}}
\newcommand{\ch}{{\mathcal{H}}}

\newcommand{\cm}{{\mathcal{M}}}
\newcommand{\ce}{{\mathcal{E}}}
\newcommand{\cv}{{\mathcal{V}}}

\newcommand{\cc}{{\mathcal{C}}}

\newcommand{\co}{{\mathcal{O}}}

\newcommand{\ct}{{\mathcal{T}}}
\newcommand{\cj}{{\mathcal{J}}}

\renewcommand{\a}{\alpha}
\renewcommand{\b}{\beta}
\renewcommand{\l}{\lambda}
\renewcommand{\ll}{\Lambda}
\newcommand{\s}{\sigma}

\newcommand{\p}{\varphi}

\renewcommand{\th}{\theta}
\renewcommand{\d}{\delta}

\renewcommand{\o}{\omega}

\newcommand{\g}{\gamma}
\renewcommand{\gg}{\Gamma}
\newcommand{\ep}{\varepsilon}
\newcommand{\z}{\zeta}

\newcommand{\ka}{\kappa}

\newcommand{\ag}{A(\gg)\,}

\newcommand{\whe}{\widehat e}

\DeclareMathOperator{\lc}{span}
\DeclareMathOperator{\ran}{rank}

\DeclareMathOperator{\sgn}{sgn}

\allowdisplaybreaks
\numberwithin{equation}{section}

\newtheorem{theorem}{Theorem}[section]

\newtheorem{corollary}[theorem]{Corollary}
\newtheorem{proposition}[theorem]{Proposition}

\theoremstyle{definition}
\newtheorem{definition}[theorem]{Definition}
\newtheorem{remark}[theorem]{Remark}

\newtheorem{example}[theorem]{Example}

\begin{document}

\title[Spectra of infinite graphs]
{Spectra of infinite graphs: two methods of computation}
\author[L. Golinskii]{L. Golinskii}

\address{B. Verkin Institute for Low Temperature Physics and
Engineering, 47 Science ave., Kharkov 61103, Ukraine}
\email{golinskii@ilt.kharkov.ua}

\date{\today}

\subjclass[2010]{}


\begin{abstract}
Two method for computation of the spectra of certain infinite graphs are suggested. The first one can be viewed as
a reversed Gram--Schmidt orthogonalization procedure. It relies heavily on the spectral theory of Jacobi matrices.
The second method is related to the Schur complement for block matrices.
A number of examples including finite graphs with tails, chains of cycles and ladders are worked out in detail.
\end{abstract}

\maketitle

\tableofcontents

\section{Introduction}
\label{s1}

We begin with some rudiments of the graph theory. For the sake of simplicity we restrict ourselves
with simple, connected, undirected, finite or infinite (countable) weighted graphs, although
the main result holds for weighted multigraphs and graphs with loops as well.
We will primarily label the vertex set $\cv(\gg)$ by positive integers $\bn=\{1,2,\ldots\}$,
$\{v\}_{v\in \cv}=\{j\}_{j=1}^\o$, $\o\le\infty$. The symbol $i\sim j$ means that the vertices
$i$ and $j$ are incident, i.e., $\{i,j\}$ belongs to the edge set $\ce(\gg)$.
A graph $\gg$ is weighted if a positive number $d_{ij}$ (weight) is assigned to each edge
$\{i,j\}\in\ce(\gg)$. When $d_{ij}=1$ for all $i,j$, the graph is unweighted.

The degree (valency) of a vertex $v\in\cv(\gg)$ is a number $\g(v)$ of edges emanating from $v$.
A graph $\gg$ is said to be locally finite, if $\g(v)<\infty$ for all $v\in\cv(\gg)$, and
uniformly locally finite, if $\sup_{\cv}\g(v)<\infty$. The latter will be the case in all our
considerations below.

The spectral graph theory deals with the study of spectra and spectral properties of certain
matrices related to graphs (more precisely, operators generated by such matrices in the standard
basis $\{e_k\}_{k\in\bn}$ and acting in the corresponding Hilbert spaces $\bc^n$,
$\ell^2=\ell^2(\bn)$, or, more generally, $\ell^2(\cv(\gg))$.
One of the most notable of them is the {\it adjacency matrix} $A(\gg)$
\begin{equation}\label{adjmat}
A(\gg)=\|a_{ij}\|_{ij=1}^\o, \quad
a_{ij}=\left\{
  \begin{array}{ll}
    d_{ij}, & \{i,j\}\in\ce(\gg); \\
    0, & \hbox{otherwise.}
  \end{array}
\right.
\end{equation}
The corresponding adjacency operator will be denoted by the same symbol. It acts as
\begin{equation}\label{adjop}
A(\gg)\,e_k=\sum_{j\sim k} a_{jk}\,e_j, \qquad k\in\bn.
\end{equation}
Clearly, $\ag$ is a symmetric, densely-defined linear operator, whose domain is the set of all finite
linear combinations of the basis vectors. The operator $A(\gg)$ is bounded and selfadjoint in $\ell^2$,
as long as the graph $\gg$ is uniformly locally finite.

Under the {\it spectrum} $\s(\gg)$ ({\it resolvent set} $\rho(\gg)$) of the graph we always mean the spectrum (the resolvent set) of its
adjacency operator $A(\gg)$. We stick to the following classification for the parts of the spectrum (one of the most confusing notions
in the spectral theory). Under the {\it discrete spectrum} $\s_d(\gg)$ we mean the set of all isolated eigenvalues of $\ag$ of finite
multiplicity. The {\it essential spectrum} $\s_{ess}(\gg)$ is the complement $\s(\gg)\backslash\s_d(\gg)$. If $H_\gg\{(a,b)\}$ is the
spectral subspace of $\ag$ for the interval $(a,b)$, a point $\l\in\s_{ess}(\gg)$ if and only if for each $\ep>0$ the dimension
$\dim H_\gg\{(\l-\ep, \l+\ep)\}=+\infty$. We denote by $\s_p(\gg)$ the set of {\it all} eigenvalues of $\ag$, called the {\it point spectrum} of $\gg$.
Sometimes we use a notation $\s_h(\gg)$ for the set of eigenvalues of $\ag$, lying on the essential spectrum (the hidden spectrum).
We will observe such phenomenon, $\s_h(\gg)\not=\emptyset$, in a number of subsequent examples.

The underlying Hilbert space, wherein the adjacency operator $\ag$ acts, is $\ell^2$ as soon as we use the set $\bn$ to enumerate the vertex set
$\cv(\gg)$. But sometimes it is much more convenient to use another set of indices, see, e.g., Example \ref{lantern}. In general, the underlying Hilbert
space $\ell^2(\gg)$ is the set of all square summable sequences defined on $\cv(\gg)$. The standard basis in this space is $\{e_v(\cdot)\}_{v\in\cv(\gg)}$ with
\begin{equation*}
e_v(w)=\left\{
            \begin{array}{ll}
              1, & w=v; \\
              0, & w\not=v.
            \end{array}
          \right.
\end{equation*}
Relation \eqref{adjop} (for unweighted graphs) looks as
\begin{equation}\label{adjop1}
(A(\gg)\,e_v)(w)=\sum_{u\sim v} e_u(w), \qquad u,v\in\cv(\gg).
\end{equation}

Whereas the spectral theory of finite graphs is very well established (see, e.g., \cite{Bap, BrHae, Chung97, CDS80}), the corresponding theory
for infinite graphs is in its infancy. We refer to \cite{M82, MoWo89, SiSz} for the basics of this theory. In contrast to the general
consideration in \cite{MoWo89}, the goal of these notes is to carry out a complete spectral analysis for certain classes of infinite graphs.

We suggest two methods of computation of spectra for such graphs. The first one applies to the graphs which can be called ``finite graphs with
tails attached to them'' and some other closely related graphs. It is pursued in two stages. At the first one, we construct a canonical model
for the adjacency operators of such graphs, which is an orthogonal sum of a finite dimensional operator and a Jacobi operator of finite rank (finite
and Jacobi components of the graph). At the second stage, the spectrum of the Jacobi component is computed by means of the Jost solution
for the corresponding recurrence relation, and the spectrum of the finite component by the standard means of linear algebra.

To be precise, we define first an operation of coupling well known for finite graphs (see, e.g., \cite[Theorem 2.12]{CDS80}).

\begin{definition}\label{coupl}
Let $\gg_k$, $k=1,2$, be two weighted graphs with no common vertices, with the vertex sets and edge
sets $\cv(\gg_k)$ and $\ce(\gg_k)$, respectively, and let $v_k\in \cv(\gg_k)$. A weighted graph
$\gg=\gg_1+\gg_2$ will be called a {\it coupling by means of the bridge $\{v_1,v_2\}$ of weight} $d$ if
\begin{equation}\label{defcoup}
\cv(\gg)=\cv(\gg_1)\cup \cv(\gg_2), \qquad \ce(\gg)=\ce(\gg_1)\cup \ce(\gg_2)\cup \{v_1,v_2\}.
\end{equation}
So, we join $\gg_2$ to $\gg_1$ by a new edge of weight $d$ between $v_2$ and $v_1$.
\end{definition}

If the graph $\gg_1$ is finite, $V(\gg_1)=\{1,2,\ldots,n\}$, and $V(\gg_2)=\{j\}_{j=n+1}^\o$,
we can with no loss of generality put $v_1=n$, $v_2=n+1$, so the adjacency matrix $A(\gg)$ can
be written as a block matrix
\begin{equation}\label{adjcoup}
A(\gg)=\begin{bmatrix}
A(\gg_1) & E_d \\
E_d^*& A(\gg_2)
\end{bmatrix}, \qquad
E_d=\begin{bmatrix}
0 & 0 & 0 & \ldots \\
\vdots & \vdots & \vdots &  \\
0 & 0 & 0 & \ldots \\
d & 0 & 0 & \ldots
\end{bmatrix},
\end{equation}
the matrix with $n$ rows and one nonzero entry. If $\gg_2=\bp_\infty\{a_j\}$, the one-sided weighted infinite path,
$a_j=d_{j,j+1}$, we can view the coupling $\gg=\gg_1+\bp_\infty\{a_j\}$ as a finite graph with the tail.
This is the class of graphs we will primarily be dealing with here.

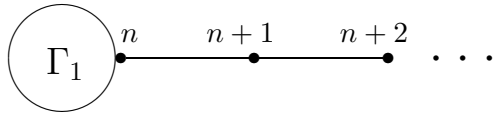
\begin{picture}(300, 100)

\put(50, 50){\circle{80}} \put(44, 44){\Large{$\Gamma_1$}}

\multiput(72, 50) (50,0) {3} {\circle* {4}}
\multiput(74, 50) (50,0) {2} {\line(1,0) {46}}
\put(72, 56) {$n$} \put(104, 56) {$n+1$} \put(154, 56) {$n+2$}
\multiput(190, 50) (10, 0) {3} {\circle*{2}}

\end{picture}

A special class of infinite matrices will play a crucial role in what follows.

Under {\it Jacobi} or {\it tridiagonal matrices} we mean here semi-infinite matrices of the form
\begin{equation}\label{defjac}
J=J(\{b_j\}, \{a_j\})_{j\in\bn}=
\begin{bmatrix}
 b_1 & a_1 & & \\
 a_1 & b_2 & a_2 & \\
     & a_2 & b_3 & \ddots & \\
     &     & \ddots & \ddots
\end{bmatrix}, \quad b_j\in\br, \quad a_j>0.
\end{equation}
They generate linear operators (called the Jacobi operators) on the Hilbert space $\ell^2(\bn)$. The matrix
\begin{equation}\label{free}
J_0:=
\begin{bmatrix}
 0 & 1 & 0 & 0 & \\
 1 & 0 & 1 & 0 & \\
 0 & 1 & 0 & 1 & \\
     & \ddots  & \ddots & \ddots & \ddots
\end{bmatrix}
\end{equation}
called a {\it discrete Laplacian} or a {\it free Jacobi matrix}, is of particular interest in the sequel. We denote by $J_{0,p}$ a $p\times p$-matrix,
which is a principal block of order $p$ of $J_0$ (the discrete Laplacian of order $p$).

The Jacobi matrices arise in the spectral graph theory because of the relation for the adjacency matrix $A(\bp_\infty\{a_j\})$ of the weighted path
\begin{equation}\label{adjpathw}
A(\bp_\infty(\{a_j\})):=J(\{0\}, \{a_j\}).
\end{equation}
In case of the unweighted path, that is, $a_j\equiv 1$, we have
\begin{equation}\label{adjpath}
A(\bp_\infty)=J_0.
\end{equation}
The spectrum of $J_0$ is $\s(J_0)=[-2,2]$.

Similarly, the discrete Laplacian of order $p$ is the adjacency matrix of the path $\bp_p$ with $p$ vertices, $J_{0,p}=A(\bp_p)$.
It is well known \cite[Section 1.4.4]{BrHae} that the spectrum
\begin{equation}\label{spfinlap}
\s(J_{0,p})=\Bigl\{2\cos\frac{\pi j}{p+1}\Bigr\}_{j=1}^p.
\end{equation}

Sometimes in our consideration two-sided Jacobi matrices $J=J(\{b_j\}, \{a_j\})_{j\in\bz}$, acting on the Hilbert space $\ell^2(\bz)$,
show up. The discrete Laplacian is
\begin{equation}\label{twosidlap}
J=J_0(\bz)=J(\{0\}, \{1\})_{j\in\bz}, \quad \s(J_0(\bz))=[-2,2].
\end{equation}

It follows from \eqref{adjcoup} that for an arbitrary finite weighted graph $G$
\begin{equation}\label{blmat1}
A(G+\bp_\infty\{a_j\})=
\begin{bmatrix}
A(G) & E_d \\
E_d^*& J(\{0\},\{a_j\})
\end{bmatrix}.
\end{equation}

To proceed further, let us recall the notions of truncation and extension for Jacobi matrices.

Given two Jacobi matrices $J_k=J(\{\b_j^{(k)}\}, \{\a_j^{(k)}\})$, $k=1,2$, the matrix $J_2$
is called a {\it truncation} of $J_1$ (and $J_1$ is an {\it extension} of $J_2$) if
\begin{equation*}
\b_j^{(2)}=\b_{j+q}^{(1)}, \quad \a_j^{(2)}=\a_{j+q}^{(1)}, \qquad j\in\bn,
\end{equation*}
for some $q\in\bn$. In other words, $J_2$ is obtained from $J_1$ by deleting the first $q$
rows and columns. The term {\it $q$-stripped matrix} is also in common usage.
If $J_2=J_0$, $J_1$ is said to be a {\it Jacobi matrix of finite rank} or an {\it eventually free Jacobi matrix}.

For Jacobi matrices $J$ of finite rank, it is well known that
\begin{equation}
\s(J)=\s_{ess}(J)\cup\s_d(J)=[-2,2]\cup\s_d(J),
\end{equation}
the discrete spectrum $\s_d(J)$ is finite, and the union is disjoint.

We suggest a ``canonical'' form for the block matrices \eqref{blmat1} and the algorithm of their reduction to this form.

\begin{theorem}\label{algorithm}
Let $A$ be a block matrix on $\ell^2$,
\begin{equation}\label{adjmat}
A=\begin{bmatrix}
\mathcal{A} & E_d \\
E_d^*& J
\end{bmatrix}, \qquad
E_d=\begin{bmatrix}
0 & 0 & 0 & \ldots \\
\vdots & \vdots & \vdots &  \\
0 & 0 & 0 & \ldots \\
d & 0 & 0 & \ldots
\end{bmatrix}, \quad d>0,
\end{equation}
where $\mathcal{A}=[a_{ij}]_{i,j=1}^n$ is a real symmetric matrix of order $n$, $J=J(\{\beta_j\}, \{\alpha_j\})$ a Jacobi matrix.
Then $A$ can be reduced to the block diagonal form
\begin{equation}\label{canform}
A\simeq \begin{bmatrix}
\widehat{\mathcal A} &  \\
  & \widehat J
\end{bmatrix},
\end{equation}
where $\widehat{\mathcal{A}}$ is a real symmetric matrix of order at most $n-2$, and the Jacobi matrix $\widehat J$ is an extension of $J$.
In other words, there is a unitary operator $U$ on $\ell^2$ such that
\begin{equation}
U^{-1}AU=\widehat{\mathcal{A}}\oplus \widehat J.
\end{equation}
\end{theorem}

\begin{corollary}
Given a finite weighted graph $G$, the adjacency operator of the coupling $\gg=G+\bp_\infty\{a_j\}$
is unitarily equivalent to the orthogonal sum
\begin{equation}
U^{-1}A(\gg)U= F(\gg)\oplus J(\gg)
\end{equation}
of a finite-dimensional operator $F(\gg)$ and a Jacobi operator $J(\gg)$, which is an extension of
$J(\{0\},\{a_j\})$.
\end{corollary}

The matrix $J(\gg)$ is of finite rank, as long as $\bp_\infty$ is unweighted. We call $F(\gg)$ a {\it finite-dimensional component} of
the coupling $\gg$, and $J(\gg)$ its {\it Jacobi component}.

The proof of Theorem \ref{algorithm} in \cite{Gol16} can be viewed as a ``reversed Gram--Schmidt algorithm''. Denote by $\{e_j\}_{j\ge1}$
the standard orthonormal basis in $\ell^2$. We construct a new orthonormal basis (canonical basis) $\{h_j\}_{j\ge1}$ so that $h_j=e_j$,
$j=n,n+1,\ldots$, and for some $q=1,\ldots,n-1$ the span of $\{h_j\}_{j\ge q}$, denoted by $\cj(\gg)$, is invariant for the adjacency operator
$A(\gg)$. What is more to the point, the restriction of $A(\gg)$ on $\cj(\gg)$ has a Jacobi matrix in the basis $\{h_j\}_{j\ge q}$.
The orthogonal complement $\cf(\gg)=\ell^2\ominus\cj(\gg)=\lc\{h_j\}_{j=1}^{q-1}$ is finite-dimensional and invariant for $A(\gg)$.
The restriction of $A(\gg)$ on $\cj(\gg)$ provides the Jacobi component, and on $\cf(\gg)$ the finite-dimensional component. The case
$q=1$ means that the finite-dimensional component is missing.

It follows from the above canonical form that the spectrum of $\gg$ is
\begin{equation*}
\s(\gg)=\s(F(\gg))\bigcup\s(J(\gg)).
\end{equation*}
Hence, to compute the spectrum of $\gg$, we apply the spectral result of Damanik and Simon for Jacobi matrices of finite rank, based on
the Jost solution. The eigenvalues of the finite-dimensional component are the roots of the corresponding characteristic polynomial.

\begin{remark}\label{multicoupl}
Given a finite graph $G$, one can attach $p\ge1$ copies of the infinite path $\bp_\infty$ to {\it some} vertex
$v\in\cv(G)$. Although the graph $\gg$ thus obtained is not exactly the coupling in the sense of
Definition \ref{coupl}, its adjacency operator acts similarly to one for the coupling. Indeed, it is not
hard to see that
\begin{equation}
A(\gg)=\begin{bmatrix}
A(G) & E_{d} \\
E_{d}^*& J_0 &
\end{bmatrix} \bigoplus\Bigl(\bigoplus_{i=1}^{p-1} J_0\Bigr), \quad d:=\sqrt{p}.
\end{equation}
Hence, Theorem \ref{algorithm} applies, and the spectral analysis of such graph can be carried out.
\end{remark}

Surprisingly enough, the case, when $p\ge1$ infinite rays are attached to {\it each} vertex of a
finite graph $G$, is easy to work out, and the spectrum of such graph can be found explicitly in terms
of the spectrum of $G$. Denote such graph by $G_\infty(p)$.

\begin{theorem}\label{sungraph}
Given a finite graph $G$ of order $n$ with $\s(G)=\{\l_j\}_{j=1}^n$, let $\gg=G_\infty(p)$, $p\in\bn$.
Denote by $J(\l_j,\sqrt{p})$ the Jacobi matrices of rank $1$
\begin{equation}
J(\l_j,\sqrt{p}):=J(\{\l_j,0,0,\ldots\},\{\sqrt{p},1,1,\ldots\}).
\end{equation}
 Then the adjacency operator $A(\gg)$ is unitarily equivalent to the orthogonal sum
\begin{equation}\label{canfor12}
A(\gg) \simeq \bigoplus_{j=1}^n J(\l_j,\sqrt{p})
\bigoplus\Bigl(\bigoplus_{i=1}^{(p-1)n}J_0\Bigr).
\end{equation}
The spectrum of $\gg$ is
\begin{equation}\label{specgen}
\s(\gg)=[-2,2]\bigcup\left(\bigcup_{j=1}^n \s_d\bigl(J(\l_j,\sqrt{p})\bigr)\right).
\end{equation}
\end{theorem}
For the proof see \cite[Theorem 1.6]{Gol16}.

\smallskip

The spectral theory of infinite graphs with one or several rays attached to certain finite graphs
was initiated in \cite{Le-star, LeNi-dan, LeNi-umzh, Niz14} wherein several particular examples
of unweighted (background) graphs are examined. The spectral analysis of similar graphs appeared earlier
in the study of thermodynamical states on complex networks \cite{FiGuIs}. We argue in the spirit of \cite{Br07, Br071, Si96} and
supplement to the list of examples. The general canonical form for the adjacency matrices of such graphs
and the algorithm of their reducing to this form suggested in the paper apply to a wide class of couplings
(not only the graphs with tails), and also to Laplacians on graphs of such type.

\smallskip

The second method, which can be called the ``Schur complement method'', is based on this well-known notion from the algebra
of block matrices. The method is applied to wider classes of infinite graphs, as well as to some other operators (Laplacians) on graphs.

Let
\begin{equation}\label{blockoper}
\ca=\begin{bmatrix}
A_{11} & A_{12} \\
A_{21} & A_{22}
\end{bmatrix}
\end{equation}
be a block operator matrix which acts on the orthogonal sum $\ch_1\oplus\ch_2$ of two Hilbert spaces.
If $A_{11}$ is invertible, the matrix $\ca$ can be factorized as
\begin{equation}\label{factor1}
\ca=\begin{bmatrix}
I & 0 \\
A_{21}A_{11}^{-1} & I
\end{bmatrix}
\begin{bmatrix}
A_{11} & 0 \\
0 & C_{22}
\end{bmatrix}
\begin{bmatrix}
I & A_{11}^{-1}A_{12} \\
0 & I
\end{bmatrix},
\end{equation}
$I$ is the unity operator on the corresponding Hilbert space. Similarly, if $A_{22}$ is invertible, one can write
\begin{equation}\label{factor2}
\ca=\begin{bmatrix}
I & A_{12}A_{22}^{-1} \\
0 & I
\end{bmatrix}
\begin{bmatrix}
C_{11} & 0 \\
0 & A_{22}
\end{bmatrix}
\begin{bmatrix}
I & 0 \\
A_{22}^{-1}A_{21} & I
\end{bmatrix}.
\end{equation}
Here
\begin{equation}\label{schurcom}
C_{22}:=A_{22}-A_{21}A_{11}^{-1}A_{12}, \qquad C_{11}:=A_{11}-A_{12}A_{22}^{-1}A_{21}
\end{equation}
are usually referred to as the {\it Schur complements} \cite{Schurcomp}, \cite[Section 0.8.5]{HoJo}.
Both equalities can be checked by direct multiplication.

The result below follows immediately from the formulae \eqref{factor1} and \eqref{factor2}.

\begin{proposition}\label{pro1}
Given a block operator matrix $\ca$ $\eqref{blockoper}$, let $A_{22}$ $(A_{11})$ be invertible.
Then $\ca$ is invertible if and only if so is $C_{11}$ $(C_{22})$.
\end{proposition}

Note that in the premises of Proposition \ref{pro1} the inverse $\ca^{-1}$ takes the form
\begin{equation*}
\ca^{-1}=\begin{bmatrix}
C_{11}^{-1} & -C_{11}^{-1}A_{12}A_{22}^{-1} \\
-A_{22}^{-1}A_{21}C_{11}^{-1} & A_{22}^{-1}+A_{22}^{-1}A_{21}C_{11}^{-1}A_{12}A_{22}^{-1}
\end{bmatrix}
\end{equation*}
and, respectively,
\begin{equation*}
\ca^{-1}=\begin{bmatrix}
A_{11}^{-1}+A_{11}^{-1}A_{12}C_{22}^{-1}A_{21}A_{11}^{-1} & -A_{11}^{-1}A_{12}C_{22}^{-1} \\
-C_{22}^{-1}A_{21}A_{11}^{-1} & C_{22}^{-1}
\end{bmatrix}.
\end{equation*}

Denote by $\rho(T)$ the resolvent set of a bounded, linear operator $T$, i.e., the set of
complex numbers $\l$ so that $\l I-T$ is boundedly invertible. We apply the latter result to the block matrix
\begin{equation}\label{blokres}
\l I-\ca=\begin{bmatrix}
\l I-A_{11} & -A_{12} \\
-A_{21} & \l I-A_{22}
\end{bmatrix}, \qquad \l\in\bc,
\end{equation}
to obtain

\begin{proposition}\label{pro2}
Given a block operator matrix $\ca$ $\eqref{blockoper}$, let $\l\in\rho(A_{22})$
$\bigl(\l\in\rho(A_{11})\bigr)$.
Then $\l\in\rho(\ca)$ if and only if the operator
\begin{equation}\label{schucom}
\begin{split}
C_{11}(\l) &=\l I-A_{11}-A_{12}(\l I-A_{22})^{-1}\,A_{21} \\
(C_{22}(\l) &=\l I-A_{22}-A_{21}(\l I-A_{11})^{-1}\,A_{12})
\end{split}
\end{equation}
is invertible.
\end{proposition}

We proceed as follows. The basics of the spectral theory for certain classes of Jacobi matrices are presented in Section \ref{s2}. 
In the next two sections we collect a number of illuminating examples, the graphs with infinite tails (Section \ref{s3}) and ladders 
and chains of cycles (Section \ref{s4}). We construct explicitly the canonical bases and find the spectra of the corresponding graphs. 
In Section \ref{s5} we discuss another method based on the Schur complement.

\section{Spectral analysis for classes of Jacobi matrices}
\label{s2}

\subsection{Jacobi matrices of finite rank and Jost functions}

For the Jacobi matrices of finite rank a complete spectral analysis is available at
the moment, see \cite{DaSi06, Ko14}. A basic object known as the {\it perturbation determinant}
\cite{GK69} is a key ingredient of perturbation theory.

Given bounded linear operators $T_0$ and $T$ on the Hilbert space such that $T-T_0$ is a nuclear
operator, the perturbation determinant is defined as
\begin{equation}\label{perdet}
L(\l;T,T_0):=\det(I+(T-T_0)R(\l,T_0)), \quad R(\l,T_0):=(T_0-\l)^{-1}
\end{equation}
is the resolvent of the operator $T_0$, an analytic operator-function on the resolvent set $\rho(T_0)$.

The perturbation determinant is designed for the spectral analysis of the perturbed operator $T$, once
the spectral analysis for the unperturbed one $T_0$ is available. In particular, the essential spectra of $T$ and $T_0$
agree, and the discrete spectrum of $T$ is exactly the zero set of the analytic function $L$ on $\rho(T_0)$,
at least if the latter is a domain, i.e., a connected, open set in the complex plane.

In the simplest case, $\ran(T-T_0)<\infty$, the perturbation determinant is the standard finite-dimensional determinant. Indeed, now
$$ (T-T_0)h=\sum_{k=1}^p \langle h, \varphi_k\rangle\,\psi_k, \quad
(T-T_0)R(\l,T_0)h=\sum_{k=1}^p \langle h, R^*(\l,T_0)\varphi_k\rangle\,\psi_k,
$$
so $L$ can be computed by the formula (see, e.g., \cite[Section IV.1.3]{GK69})
\begin{equation}\label{compd}
L(\l; T,T_0)=\det[\delta_{ij}+\langle R(\l,T_0)\psi_i,\varphi_j\rangle]_{i,j=1}^p.
\end{equation}

Our particular concern is $T_0=J_0$, the free Jacobi matrix. The matrix of its resolvent in the
standard basis in $\ell^2$ is given by (see, e.g., \cite{KiSi03})
\begin{equation}\label{resfree}
R(\l, J_0)=[r_{ij}(z)]_{i,j=1}^\infty, \ \
r_{ij}(z)=\frac{z^{|i-j|}-z^{i+j}}{z-z^{-1}}\,,  \ \ \l=z+\frac1z, \ \ z\in\bd.
\end{equation}
If $T=J$ is a Jacobi matrix of finite rank $p$, we end up with computation
of the ordinary determinant \eqref{compd} of order $p$.

It is instructive for the further usage computing two simplest perturbation determinants for
$\ran\,(J-J_0)=1$ and $2$.

\begin{example}\label{pd1}
Let
$$ J=J(\{b_j\}, \{1\}): \quad b_j=0, \quad j\not=q,
$$
so $J-J_0=\langle\cdot,e_q\rangle\,b_q e_q$. By \eqref{resfree} and \eqref{compd},
\begin{equation}\label{rank1}
\widehat L(z) :=L\Bigl(z+\frac1{z}\,; J,J_0\Bigr)=1+b_q r_{qq}(z)
=1-b_q z\,\frac{z^{2q}-1}{z^2-1}\,.
\end{equation}

Similarly, let
$$ J=J(\{0\}, \{a_j\}): \quad a_j=1, \quad j\not=q,
$$
so $J-J_0=\langle\cdot,e_q\rangle\,(a_q-1)\,e_{q+1}+\langle\cdot,e_{q+1}\rangle\,(a_q-1)\,e_{q}$,
and again
\begin{equation}\label{rank2}
\begin{split}
\widehat L(z) &=\begin{vmatrix}
1+(a_q-1)\, r_{q,q+1}(z) & (a_q-1)\, r_{qq}(z) \\
(a_q-1)\, r_{q+1,q+1}(z) & 1+(a_q-1)\, r_{q+1,q}(z)
\end{vmatrix}
\\
&=1+(1-a_q^2) z^2\,\frac{z^{2q}-1}{z^2-1}\,.
\end{split}
\end{equation}
\end{example}

\smallskip

In the Jacobi matrices setting there is yet another way of computing perturbation determinants
based on the so-called Jost solution and Jost function (see, e.g., \cite[Section 3.7]{SiSz}).

Consider the basic recurrence relation for the Jacobi matrix $J$
\begin{equation}\label{3term}
a_{n-1}y_{n-1}+b_ny_n+a_ny_{n+1}=\Bigl(z+\frac1{z}\Bigr)\,y_n,  \quad z\in\bd, \quad n\in\bn,
\end{equation}
where we put $a_0=1$. Its solution $y_n=u_n(z)$ is called the {\it Jost solution} if
\begin{equation}\label{js1}
\lim_{n\to\infty} z^{-n}u_n(z)=1, \qquad z\in\bd.
\end{equation}
In this case the function $u=u_0(z)$ is called the {\it Jost function}.

The Jost solution certainly exists for finite rank Jacobi matrices. The Jost function is
now an algebraic polynomial, called the {\it Jost polynomial}.
Indeed, let
$$ a_q\not=1, \quad a_{q+1}=a_{q+2}=\ldots=1, \qquad b_{q+1}=b_{q+2}=\ldots=0. $$
One can put $u_k(z)=z^k$, $k=q+1,q+2,\ldots$ and then determine $u_q, u_{q-1},\ldots,u_0$
consecutively from \eqref{3term}. So,
\begin{equation}\label{js2}
\begin{split}
a_q\,u_q(z) &=z^q, \\
a_{q-1}a_q\,u_{q-1}(z) &=\a_q\,z^{q+1}-b_qz^q+z^{q-1}, \quad \a_q:=1-a_q^2,
\end{split}
\end{equation}
etc., and in general, for $k=0,\ldots,q$
$$ u_{q-k}(z)=\sum_{j=-k}^k \b_{q,j}z^{q+j}, \ \ \b_{q,j}\in\br, \ \ \b_{q,-k}=1, \ \ \b_{q,k}=\frac{1-a_q^2}{a_{q-k}a_{q-k+1}\ldots a_q}\,. $$
In particular, for $q=1$
\begin{equation}\label{jf1}
a_1\,u(z)=\a_1\,z^2-b_1z+1,
\end{equation}
and for $q=2$
\begin{equation}\label{jf2}
a_1a_2\,u(z) =\a_2\,z^4-(b_2+b_1\a_2)\,z^3 +(\a_1+\a_2+b_1b_2)\,z^2 -(b_1+b_2)\,z+1.
\end{equation}

The relation between the perturbation determinant and the Jost function is given by
\begin{equation}\label{pdfj}
u(z)=\prod_{j=1}^\infty a_j^{-1}\cdot \widehat L(z),
\end{equation}
see \cite{KiSi03}, and such recursive way of computing perturbation determinants is sometimes far
easier than computing ordinary determinants \eqref{compd}, especially for large enough ranks of perturbation.
On the other hand, for small ranks of perturbation (as in Example \ref{pd1}) with large $q$ it is much easier
applying formula \eqref{compd}.

\begin{example}\label{jfcomp}
Let $J=J(\{b_j\}, \{a_j\})$ be a Jacobi matrix such that
$$ b_j=0, \quad j\not=1; \qquad a_j=1, \quad j\not=q. $$
We have $u_{q+j}(z)=z^{q+j}$, $j=1,2\ldots$,
\begin{equation*}
\begin{split}
a_q u_q(z) &=z^q, \quad a_qu_{q-1}(z) =\a_q\,z^{q+1}+z^{q-1}, \\
a_qu_{q-2} &=\a_q\,(z^{q+2}+z^q)+z^{q-2},
\end{split}
\end{equation*}
and, by the induction,
\begin{equation}\label{jf21}
a_qu_{q-k}(z)=\a_q\,z^{q-k+2}\,\frac{z^{2k}-1}{z^2-1}\,,\qquad k=1,2,\ldots,q-1.
\end{equation}

Next, for $q=1$ we have exactly \eqref{jf1}, so let $q\ge2$.
The recurrence relation \eqref{3term} with $n=1$ gives
$$ a_q u(z)+b_1a_q u_1(z)+a_q u_{2}(z)=\Bigl(z+\frac1{z}\Bigr)\,a_qu_1(z), $$
and so we come to the following expression for the Jost polynomial
\begin{equation}\label{jf3}
a_qu(z)(z^2-1)=\a_q\,(z-b_1)\,z^{2q+1}-b_1a_q^2 z^3+a_q^2 z^2+b_1 z-1.
\end{equation}

Similarly, for the Jacobi matrix $J=J(\{b_j\}, \{a_j\})$ with
$$ b_j=0, \ \  j\not=q; \quad a_j=1, \quad j\not=1 $$
one has
\begin{equation}\label{jf4}
a_1u(z)=-b_q\,\frac{z^{2q+1}+\a_1z^{2q-1}-\a_1z^3-z}{z^2-1}+\a_1z^2+1.
\end{equation}
For the Jacobi matrix $J=J(\{0\}, \{a_j\})$ with $a_j=1$, $j\not=1$, $q$,
the Jost polynomial is given by
\begin{equation}\label{jf5}
a_1a_q u(z)=\a_q\,\frac{z^{2q+2}+\a_1z^{2q}-\a_1z^4-z^2}{z^2-1}+ \a_1z^2+1.
\end{equation}
\end{example}

\medskip

The spectral theorem for finite rank Jacobi matrices due to Damanik and Simon \cite{DaSi06}
provides a complete description of the spectral measure of such matrices.

{\bf Theorem} (Damanik--Simon). Let $J=J(\{b_j\}, \{a_j\})$ be a Jacobi matrix of finite rank
$$ a_{q+1}=a_{q+2}=\ldots=1, \qquad  b_{q+1}=b_{q+2}=\ldots=0, $$
and $u=u_0(J)$ be its Jost polynomial. Then
\begin{itemize}
  \item $u$ is a real polynomial of degree $\deg u\le 2q$,
  $\deg u=2q$ if and only if $a_q\not=1$.
  \item All roots of $u$ in the unit disk $\bd$ are real and simple, $u(0)\not=0$. A number $\l_j$
  is an eigenvalue of $J$ if and only if
\begin{equation}\label{zhuk}
\l_j=z_j+\frac1{z_j}, \qquad z_j\in(-1,1), \quad u(z_j)=0.
\end{equation}
  \item The spectral measure $\s(J)$ is of the form
\begin{equation}
\s(J,dx)=\s_{ac}(J,dx)+\s_d(J,dx)=w(x)\,dx+\sum_{j=1}^N \s_j\d(\l_j),
\end{equation}
where
$$ w(x):=\frac{\sqrt{4-x^2}}{2\pi|u(e^{it})|^2}\,, \quad x=2\cos t, \quad \s_j=\frac{z_j(1-z_j^{-2})^2}{u'(z_j)u(1/z_j)}\,. $$
\end{itemize}
Note that $|u(e^{it})|^2=Q(x)$, $x=2\cos t$, $Q$ is a real polynomial of the same degree as the
Jost polynomial $u$.

\smallskip

The algebraic equations which we encounter later on cannot in general be solved explicitly.
By means of the following well-known result \cite[p.~41]{PoSz}, we can
determine how many roots (if any) they have in $(-1,1)$.

{\bf Theorem} (Descarte's rule). Let $a(x)=a_0x^n+\ldots+a_n$ be a real polynomial. Denote by $\mu(a)$
the number of its positive roots, and $\nu(a)$ the number of the sign changes in the sequence
$\{a_0,\ldots,a_n\}$ of its coefficients (the zero coefficients are not taken into account).
Then $\nu(a)-\mu(a)$ is a nonnegative even number.

\subsection{Periodic Jacobi matrices}

See \cite[Chapter 7]{Tjo}, \cite[Chapter 5]{SiSz} for an extensive theory.

Let us go back to the basic recurrence relations \eqref{3term}
\begin{equation}\label{3term1}
a_{n-1}y_{n-1}(\l)+b_ny_n(\l)+a_ny_{n+1}(\l)=\l\,y_n,  \quad n\in\bn, \quad a_0=1,
\end{equation}
and consider its two solutions
$$ y_n=p_n:\ \ p_0=0, \ p_1=1; \qquad y_n=q_n:\ \ q_0=-1, \ q_1=0. $$
The polynomials $p_n$ $(q_n)$ are called the first (second) kind polynomials for the Jacobi matrix $J$. We have
\begin{equation}\label{1kind}
p_2(\l)=\frac{\l-b_1}{a_1}\,, \quad p_3(\l)=\frac{(\l-b_2)(\l-b_1)}{a_1a_2}-\frac{a_1}{a_2}\,, \ldots,
\end{equation}
so $\deg p_k=k-1$.

Recall that $1$-stripped matrix $J_1$ for $J$ is given by $J_1=J(\{b_{n+1}\}, \{a_{n+1}\})$. The stripping formula \cite[formula (3.2.16)]{SiSz} relates
the second kind polynomials for $J$ and the first kind ones for $J_1$
\begin{equation}\label{strip}
q_n(\l)=\frac{p_{n-1}^{(1)}(\l)}{a_1}, \qquad n\in\bn,
\end{equation}
so $\deg q_k=k-2$.

\begin{example}\label{cheb} ``Chebyshev polynomials''.

We compute the 1st (2nd) kind polynomials for two particular Jacobi matrices. Recall the notion of Chebyshev polynomials of the 1st and
2nd kind
\begin{equation}\label{cheb}
T_n(\cos\th)=\cos n\th, \quad U_n(\cos\th)=\frac{\sin(n+1)\th}{\sin\th}\,, \qquad n=0,1,\ldots,
\end{equation}
so $T_0=U_0=1$,
\begin{equation*}
\begin{split}
T_1(\l) &=\l, \ \ \ \  T_2(\l)=2\l^2-1, \ \ T_3(\l)=4\l^3-3\l, \ \ T_4(\l)=8\l^4-8\l^2+1; \\
U_1(\l) &=2\l, \ \ U_2(\l)=4\l^2-1, \ \ U_3(\l)=8\l^3-4\l \ \ U_4(\l)=16\l^4-12\l^2+1.
\end{split}
\end{equation*}
The general expressions as the products are
$$ T_n(\l)=2^{n-1}\,\prod_{k=1}^n \Bigl(\l-\cos\frac{(2k-1)\pi}{2n}\Bigr), \quad
U_n(\l)=2^{n}\,\prod_{k=1}^n \Bigl(\l-\cos\frac{k\pi}{n+1}\Bigr), $$

The standard equalities
\begin{equation*}
\begin{split}
\cos(n-1)\th+\cos(n+1)\th &= 2\cos\th\,\cos n\th, \\
\sin(n-1)\th+\sin(n+1)\th &= 2\cos\th\,\sin n\th
\end{split}
\end{equation*}
lead to
\begin{equation*}
\begin{split}
T_{n-1}(\l)+T_{n+1}(\l) &=2\l\,T_n(\l), \\
U_{n-1}(\l)+U_{n+1}(\l) &=2\l\,U_n(\l), \quad n\in\bn.
\end{split}
\end{equation*}
It is clear now that the 1st kind polynomials for $J=J_0$ are
$$ p_0(J_0)=0, \quad p_n(\l,J_0)=U_{n-1}\Bigl(\frac{\l}2\Bigr), \quad n=1,2,\ldots, $$
and 1st kind polynomials for $J=J_0'=J(\{0\}, \{\sqrt2,1,1,\ldots\})$ are
$$ p_0(J_0')=0, \ \ p_1(\l,J_0')=1, \ \  p_n(\l,J_0')=\sqrt2\,T_{n-1}\Bigl(\frac{\l}2\Bigr), \quad n=2,3,\ldots. $$

More generally, if $b_1=b_2=\ldots=0$, $a_1=a_2=\ldots=a_k=1$, then
\begin{equation}\label{1stkind1}
p_n(\l,J)=U_{n-1}\Bigl(\frac{\l}2\Bigr), \quad n=1,2,\ldots, k+1.
\end{equation}
If $b_1=b_2=\ldots=0$, $a_1=\sqrt2$, $a_2=\ldots=a_k=1$, then
\begin{equation}\label{1stkind2}
p_n(\l,J)=\sqrt2\,T_{n-1}\Bigl(\frac{\l}2\Bigr), \quad n=2,3,\ldots, k+1. \end{equation}
\end{example}

\smallskip

The matrix form of \eqref{3term1} looks
$$
\begin{bmatrix}
y_{n+1} \\
a_ny_n
\end{bmatrix}=A(\l; a_n,b_n)
\begin{bmatrix}
y_{n} \\
a_{n-1}y_{n-1}
\end{bmatrix}, \quad
A(\l; a_n,b_n)=\begin{bmatrix}
\frac{\l-b_n}{a_n} & -\frac1{a_n} \\
a_n & 0
\end{bmatrix}, \ \ n\in\bn.
$$
The composition of the latter equalities leads to the transfer matrix
\begin{equation}\label{trans}
\ct_n(\l):=A(\l; a_n,b_n)A(\l; a_{n-1},b_{n-1})\ldots A(\l; a_1,b_1).
\end{equation}
Precisely,
\begin{equation}\label{trans1}
\begin{split}
\begin{bmatrix}
p_{n+1} \\
a_np_n
\end{bmatrix} &=\ct_n(\l)
\begin{bmatrix}
1 \\
0
\end{bmatrix}, \quad
\begin{bmatrix}
q_{n+1} \\
a_nq_n
\end{bmatrix} =\ct_n(\l)
\begin{bmatrix}
0 \\
-1
\end{bmatrix}, \\
\ct_n(\l) &=\begin{bmatrix}
p_{n+1}(\l) & -q_{n+1}(\l) \\
a_np_n(\l) & -a_nq_n(\l)
\end{bmatrix}, \ \ n\in\bn.
\end{split}
\end{equation}

Since $\det A(\l; a_n,b_n)=\det \ct_n(\l)=1$, we see that for each $n\in\bn$ and complex $\l$
\begin{equation}\label{wrons}
a_n\bigl(p_n(\l)q_{n+1}(\l)-p_{n+1}(\l)q_n(\l)\bigr)=1.
\end{equation}

\smallskip

Let $\mu$ be the spectral measure of the Jacobi operator $J$. The {\it Weyl function} is defined by
$$ m(\l)=m(\l,J):=\int_\br \frac{\mu(dt)}{t-\l}\,. $$
For the initial data
$$ v_0(\l)=\begin{bmatrix}
m(\l) \\ -1 \end{bmatrix} $$
the solution of the matrix form of \eqref{3term1} is
$$ v_n(\l)=\ct_n(\l)v_0(\l)=\begin{bmatrix}
m(\l)p_{n+1}(\l)+q_{n+1}(\l) \\
a_n\bigl(m(\l)p_{n+1}(\l)+q_{n+1}(\l)\bigl)
\end{bmatrix}. $$
As is known, this solution is square summable for each $\l\in\rho(J)$. Denote $w_n:=mp_{n+1}+q_{n+1}$.
The Weyl function for the $k$-stripped Jacobi matrix $J_k$ can be expressed in terms of $w_k$ as
\begin{equation}\label{weylstr}
m(\l,J_k)=-\frac{w_k(\l)}{a_kw_{k-1}(\l)}=-\frac{m(\l)p_{k+1}(\l)+q_{k+1}(\l)}{a_k\bigl(m(\l)p_{k}(\l)+q_{k}(\l)\bigl)}\,.
\end{equation}
In particular, for $k=1$ we have
$$ m(\l,J)=\frac1{b_1-\l-a_1^2\,m(\l,J_1)}\,. $$

\smallskip

A Jacobi matrix $J$ is called {\it $N$-periodic}, $N\in\bn$, if
\begin{equation*}
a_{n+N}=a_n, \quad b_{n+N}=b_n, \qquad n=1,2,\ldots.
\end{equation*}
In other word, $J$ is $N$-periodic if and only if the $N$ stripped matrix $J_N=J$. Equality \eqref{weylstr} shows that the Weyl function satisfies
the following quadratic equation
\begin{equation}\label{quadeq}
a_Np_N(\l)\,m^2(\l)+\bigl(p_{N+1}(\l)+a_Nq_N(\l)\bigr)\,m(\l)+q_{N+1}(\l)=0.
\end{equation}
The discriminant of this equation equals in view of \eqref{wrons}
\begin{equation}\label{discrquad}
\begin{split}
D(\l) &=\bigl(p_{N+1}(\l)+a_Nq_N(\l)\bigr)^2-4a_N\,p_N(\l)q_{N+1}(\l) \\
&=\bigl(p_{N+1}(\l)-a_Nq_N(\l)\bigr)^2-4=\cd^2(\l)-4,
\end{split}
\end{equation}
where the polynomial
\begin{equation}\label{discr}
\cd(\l):=p_{N+1}(\l)-a_Nq_N(\l)=p_{N+1}(\l)-\frac{a_N}{a_1}\,p_{N-1}^{(1)}(\l)
\end{equation}
is the well-known {\it discriminant} of the Jacobi matrix $J$, which plays a key role in the spectral theory of periodic Jacobi matrices.

The ``right'' root of \eqref{quadeq} for the Weyl function can be singled out from the condition $m(z)=O(z^{-1})$, $z\to\infty$,
\begin{equation}\label{weyl}
m(\l)=\frac{-\g_N(\l)+\sqrt{\cd^2(\l)-4}}{2a_N\,p_N(\l)}\,, \ \ \g_N(\l):=p_{N+1}(\l)+\frac{a_N}{a_1}\,p_{N-1}^{(1)}(\l).
\end{equation}

\smallskip

Let us now turn to the structure of the spectrum of a periodic Jacobi matrix. It is well known, that the essential spectrum of a periodic
Jacobi matrix is of banded structure, i.e., it is a union of nondegenerate closed intervals, some of them may touch each other. Precisely, let
$J$ be an $N$-periodic Jacobi matrix. Then \cite[Section 5.4]{SiSz}
\begin{equation}\label{spess}
\begin{split}
\s_{ess}(J) &=\{x\in\br:\ -2\le\cd(x)\le2\}=\bigcup_{j=1}^N [\a_j, \b_j], \\
\a_1 &<\b_1\le\a_2<\b_2\le\ldots\le\a_N<\b_N.
\end{split}
\end{equation}
The open intervals $(\b_j, \a_{j+1})$, $j=1,2,\ldots,N-1$ are called the {\it spectral gaps}. The gap $(\b_j, \a_{j+1})$ is open
as long as $\b_j<\a_{j+1}$, and it is closed otherwise (it is convenient to count closed gaps as the ``real ones''). We have
\begin{equation*}
\begin{split}
\cd(\l) &=2\ \Leftrightarrow \ \ \ \l=\b_N, \a_{N-1}, \b_{N-2},\ldots \\
\cd(\l) &=-2\ \Leftrightarrow \ \l=\a_N, \b_{N-1}, \a_{N-2},\ldots.
\end{split}
\end{equation*}

The rest of the spectrum of $J$ is the discrete spectrum, which consists of a finite number of eigenvalues. The eigenvalues of $J$ agree with the poles of the
Weyl function \eqref{weyl}. So, the discrete spectrum $\s_d(J)$ is a part of the zero set of the polynomial $p_N$ (the denominator of \eqref{weyl}).
It is known, that the closure of each spectral gap (including the closed ones) contains exactly one root of $p_N$. To specify those, which produce the
eigenvalues, we should first choose the roots {\it inside} the gaps. Next, if $p_N(\l_0)=0$, and $\l_0$ lies inside the gap, we need the numerator in \eqref{weyl}
be {\it nonzero}, so not to cancel the root of the denominator. By \eqref{discrquad},
$$ \cd^2(\l_0)-4=\g_N^2(\l_0) \ \Rightarrow \ |\sqrt{\cd^2(\l_0)-4}|=|\g_N(\l_0)|, $$
so $\l_0$ is the eigenvalue of $J$ if and only if
\begin{equation}\label{signs}
\sgn\g_N(\l_0)=-\sgn\sqrt{\cd^2(\l_0)-4}.
\end{equation}
To determine the right side in \eqref{signs}, we enumerate all the gaps (including the closed ones), from the right to the left, with the numbers $1,2,\ldots,N-1$,
and assign the sign $-1$ to the odd numbered gaps, and the sign $1$ to the even numbered ones.

So, to find the spectrum $\s(J)=\s_{ess}(J)\cup\s_d(J)$, we proceed in three steps.

\noindent
Step 1. Find the 1st kind polynomials $p_{N+1}$ for $J$ and $p_{N-1}^{(1)}$ for $1$-stripped matrix $J_1$. Compute the discriminant $\cd(J)$ in \eqref{discr} and
the polynomial $\g_N$ in \eqref{weyl}.

\noindent
Step 2. Find the essential spectrum by \eqref{spess}.

\noindent
Step 3. Find all roots of the polynomial $p_N$ and choose those inside the gaps, for which \eqref{signs} holds.

\subsection{Right limits and eigenvalues}

Our first topic here concerns the notion of a right limit, see \cite[Chapter 7]{SiSz}.

Let $f=\{f_n\}_{n\ge1}\in\ell^\infty(\bn)$. A two-sided sequence $\p=\{\p_n\}_{n\in\bz}$ is said to be the {\it right limit} for $f$, $\p\in RL(f)$, if there is
a sequence if indices $\cm=\{m_j\}_{j\ge1}\subset\bn$ so that
\begin{equation}\label{rightlim}
\p_n=\lim_{j\to\infty} f_{n+m_j} \quad \forall n\in\bz.
\end{equation}
Sometimes we write \eqref{rightlim} as $\p_n=\lim_{m\in\cm} f_{n+m}$. Note that, although the individual value $f_{n+m_j}$ may be senseless for ``large enough''
negative $n$, the $\p_n$ in \eqref{rightlim} is well-defined. We say that $\cm$ generates the right limit $\p$.

A simple compactness argument implies the following result.

\begin{proposition}\label{prright1}
For an arbitrary $f\in\ell^\infty(\bn)$, the set $RL(f)$ is nonempty. Moreover, for each sequence of indices $\ll\subset\bn$
there is a subsequence $\cm\subset\ll$ so that \eqref{rightlim} holds.
\end{proposition}

It is clear that the set $RL(f)$ is closed under the shift $S$
\begin{equation}\label{rightshif}
\p\in RL(f) \ \Leftrightarrow \ S^k\p=\{\p_{n+k}\}_{n\in\bz}\in RL(f) \quad \forall k\in\bz.
\end{equation}

It follows directly from the definition, that
$$ \lim_{n\to\infty} f_n=g \ \Rightarrow \ RL(f)=\{\ldots, g,g,g,\ldots\}, $$
so there is one, constant sequence, in the right limit for $f$. More generally,
$$ \lim_{n\to\infty}(\tilde f_n-f_n)=0 \ \Rightarrow \ RL(\tilde f)=RL(f). $$
In particular, if $f$ and $\tilde f$ agree from some point on, the sets of right limits are the same, $RL(f)=RL(\tilde f)$.

Denote by $L(f)\subset\bc$ the set of all limit points of $f$. Clearly, $\p_n\in~L(f)$ for each $n\in\bz$ as soon as $\p\in RL(f)$. If $f$ does not converge,
the cardinality of $L(f)$ is at least $2$. It is easy to see from the second statement of Proposition~\ref{prright1} that the following holds.

\begin{proposition}\label{prright2}
Given $g_1,g_2\in L(f)$, $g_1\not=g_2$, for each $k\in\bz$ there are $\p(j)\in RL(f)$, $j=1,2$ so that
$$ \p_k(1)=g_1, \qquad \p_k(2)=g_2. $$
In particular, uniqueness of the right limit yields the convergence of $f$.
\end{proposition}

\begin{example}\label{period}
Let $f$ be an $p$-periodic sequence, $f_{n+p}=f_n$, $n\in\bn$. We extend it as $p$-periodic two-sided sequence
$$ \p_{n+p}(0)=\p_n(0), \quad n\in\bz, \qquad \p_n(0)=f_n, \quad n\in\bn, $$
so $\p(0)=\{\ldots,f_1,f_2,\ldots,f_p,f_1,f_2,\ldots\}$. It is not hard to show that
$$ RL(f)=\{S^q\p(0)\}_{q=0}^{p-1}, $$
so the right limits are exhausted by the shifts of $\p(0)$.
\end{example}

There is another, opposite in a sense, situation when the set of right limits is available.

\begin{example}\label{sparse}
An increasing sequence of positive numbers $\ll=\{\l_j\}_{j\ge1}$ is called {\it sparse}, if
$$ \lim_{i\to\infty}\l_{i+1}-\l_i=+\infty. $$
Given a sparse sequence of indices $\ll\in\bn$, put
$$ b=\{b_n\}_{n\ge1}=\left\{
                                   \begin{array}{ll}
                                     \b_1, & \hbox{$n\in\ll$;} \\
                                     \b_0, & \hbox{$n\not=\ll$.}
                                   \end{array}
                                 \right., \quad \b_0,\b_1\in\bc.
$$
We want to describe the set $RL(b)$.

Let us show first that $b^{(0)}=\{\ldots,\b_0,\b_0,\b_0,\ldots\}\in RL(b)$. Define the sequence of indices $m_j:=\l_j+(-1)^j$ for $j\ge j_1$
(and in an arbitrary way for smaller $j$). Then for each $n\in\bz$ we have $n+m_j\notin\ll$ for $j\ge j_2(n)$. Indeed, assume that
$n+m_j=\l_{s(j)}$ for an infinite number of $j$'s. Then
$$ n+(-1)^j=\l_{s(j)}-\l_j $$
for such values of $j$, that contradicts the sparseness of $\ll$.

Next, since $\b_1\in b'$, then, by Proposition \ref{prright2}, for each $k\in\bz$ there is a right limit $\p(k)\in RL(b)$ so that $\p_k(k)=\b_1$. The latter means
that there is a sequence $\{m_j\}_{j\ge1}$ (generating $\p(k)$), which satisfies
$$ b_{k+m_j}=\b_1, \quad j\ge j_3(k) \ \sim \ k+m_j\in\ll, \quad j\ge j_3(k). $$
Take $k'\not=k$ and show that $k'+m_j\notin\ll$ for $j\ge j_4(k')$. Indeed, assume on the contrary, that
$$ k'+m_j\in\ll, \quad k'+m_j=\l_{r(j)} $$
for an infinite number of $j$'s. But $k'+m_j=k'-k+k+m_j=k'-k+\l_{t(j)}$, and so
$$ k'-k=\l_{r(j)}-\l_{t(j)} $$
for such values of $j$, that contradicts sparseness of $\ll$. Hence, $\p_{k'}(k)=\b_0$ for $k'\not=k$, which means that
$$ \p(k)=b_k^{(1)}:=\{\ldots,\b_0,\b_0,\b_1,\b_0,\b_0,\ldots\}, \quad k\in\bz, $$
$\b_1$ occurs at $k$-th place. Finally,
\begin{equation}\label{rlim1}
RL(b)=\{b^{(0)}; b_k^{(1)}, \ k\in\bz\}.
\end{equation}

In exactly the same way we can examine a union of two sparse sequences. Precisely, let $R=\{r_j\}_{j\ge1}$, be a sparse sequence of positive integers,
and assume for simplicity, that $r_{i+1}-r_i$ is strictly increasing, and $r_{i+1}-r_i\ge2$. Consider the union
$$ \tilde{R}=R\cup\{R+1\}=\{r_1,r_1+1,r_2,r_2+1,\ldots\}, \quad a_n=\left\{
                                   \begin{array}{ll}
                                     \a_1, & \hbox{$n\in\tilde{R}$;} \\
                                     \a_0, & \hbox{$n\not=\tilde{R}$.}
                                   \end{array}
                                 \right.,
$$
$a=\{a_n\}_{n\ge1}$. The similar reasoning leads to the following conclusion
\begin{equation}\label{rlim2}
\begin{split}
RL(a) &=\{a^{(0)}; a_k^{(1)}, k\in\bz\}, \quad a^{(0)}=\{\ldots,\a_0,\a_0,\a_0,\ldots\},   \\
a_k^{(1)} :&=\{\ldots,\a_0,\a_0,\a_0,\a_1,\a_1,\a_0,\a_0,\a_0,\ldots\},
\end{split}
\end{equation}
$\a_1$ occurs at the places $k$, $k+1$.
\end{example}

\smallskip

Going back to Jacobi matrices $J=J(\{b_n\},\{a_n\})_{n\ge1}$, we say that a two-sided Jacobi matrix
$$J_{right}=J\bigl(\{b_n^{(r)}\},\{a_n^{(r)}\}\bigr)_{n\in\bz} $$
is a {\it right limit} of $J$ if for some sequence of indices $\{m_j\}_{j\ge1}$
$$ \lim_{j\to\infty} a_{n+m_j}=a_n^{(r)}, \quad \lim_{j\to\infty} b_{n+m_j}=b_n^{(r)}, \quad \forall n\in\bz. $$

Dealing with certain sparse graphs, we will encounter the following Jacobi matrices
\begin{equation}\label{jacsimlad}
J^{\pm}=J(\{b_n^{\pm}\}, \{1\})_{n\ge1}, \quad b_n^\pm=\left\{
                                   \begin{array}{ll}
                                     \pm1, & \hbox{$n\in\ll$;} \\
                                     0, & \hbox{$n\not=\ll$.}
                                   \end{array}
                                 \right.
\end{equation}
and
\begin{equation}
J(\{0\}, \{a_n\})_{n\ge1}, \quad a_n=\left\{
                                   \begin{array}{ll}
                                     \sqrt2, & \hbox{$n\in R\cup\{R+1\}$;} \\
                                     1, & \hbox{$n\notin R\cup\{R+1\}$.}
                                   \end{array}
                                 \right.,
\end{equation}
$\ll$ and $R$ being sparse sequences of indices. The set of right limits in the first case is given by
\begin{equation}\label{rlimsim}
RL(J^\pm)=\bigl\{J_0(\bz); J(\{\ldots,0,0,\pm1,0,0,\ldots\}, \{1\})\bigr\},
\end{equation}
$\pm1$ occurs at $k$-th place, $k\in\bz$. In the second case
\begin{equation}\label{rlimcycles}
RL(J)=\bigl\{J_0(\bz); J(\{0\}, \{\ldots,1,1,1,\sqrt2,\sqrt2,1,1,1,\ldots\})\bigr\},
\end{equation}
$\sqrt2$ occurs at the places $k$, $k+1$; $k\in\bz$.

\smallskip

The following result of Last--Simon \cite{lasi06}, \cite[Theorem 7.2.1]{SiSz}, plays a key role for computation in Section \ref{ladcycl}.

{\bf Theorem} (Last--Simon). Let $J=J(\{b_n\},\{a_n\})_{n\ge1}$ be a Jacobi matrix with bounded entries
$$ \sup_n (|a_n|+|b_n|)<\infty. $$
Then
$$ \s_{ess}(J)=\bigcup_{J_{right}\in RL(J)} \s(J_{right}). $$

\medskip

Our second topic here concerns the eigenvalues of Jacobi matrices, in particular, a result of Simon--Stolz \cite{sist96} which provides
a condition for a real $\l$ not to be an eigenvalue of $J$ ($\l\notin\s_p(J)$).
The condition is given in terms of asymptotic behavior for the norms of the transfer matrices $\ct_n$ \eqref{trans}.

Let $J=J(\{b_n\},\{a_n\})$ be a Jacobi matrix with bounded entries
$$ \sup_n (|a_n|+|b_n|)<\infty. $$
Let $\{p_n\}_{n\ge1}$ be the 1st kind polynomials for $J$. By definition, $\l\in\s_p(J)$ is equivalent to $\{p_n\}_{n\ge1}\in\ell^2$.

{\bf Theorem} (Simon--Stolz). A real number $\l\notin\s_p(J)$ as long as
\begin{equation}\label{simsto}
\sum_{n=1}^\infty \|\ct_n(\l)\|^{-2}=+\infty.
\end{equation}

The argument is simple enough. Note first, that for an invertible $2\times 2$ matrix $C$ the following holds
\begin{equation}\label{norminv}
\|C^{-1}\|=\frac{\|C\|}{|\det C|},
\end{equation}
see, e.g., \cite[Lemma 10.5.1]{Siopuc2}. Next, by \eqref{trans1} and \eqref{norminv} with $\det\ct_n=1$,
$$ 1\le \|\ct_n^{-1}(\l)\|^2\,(p_{n+1}^2(\l)+a_n^2 p_n^2(\l))\le \|\ct_n(\l)\|^2\,(1+a^2)(p_{n+1}^2(\l)+ p_n^2(\l)). $$
Hence,
$$ \frac1{1+a^2}\,\sum_{n=1}^k \|\ct_n(\l)\|^{-2}\le 2\sum_{n=1}^{k+1} p_n^2(\l), $$
and we are done.

In some examples below the condition \eqref{simsto} can be exploited.

\section{Spectra of graphs with infinite tails}
\label{s3}

To obtain the canonical form for particular graphs one should apply the reversed Gram--Schmidt algorithm by hand. Its efficiency
strongly depends on complexity of the graph in question.

\subsection{Trees with tails}

The examples below are taken partially from \cite[Section 4]{Gol16}. Some of them are new.

\begin{example}\label{starsimp} ``A weighted star''.

Let $S_n(w)$ be a simple weighted star graph of order $n+1$, $n\ge2$, with vertices $1,\ldots,n$ of degree $1$, and the vertex $n+1$ of degree $n$ being a root.
The weight of the edge $(i,n+1)$ equals $w_i$, $1\le i\le n$. We consider the coupling $\gg=S_n(w)+\bp_\infty$, where the infinite ray is
attached to the root.

The canonical basis $\{h_j\}_{j\in\bn}$ looks as follows. We put
$$ h_j:=e_j, \quad j\ge n+1 \ \Longrightarrow \ \ag h_j=h_{j-1}+h_{j+1}, \quad j\ge n+2. $$
Next, let $w:=(w_1,w_2,\ldots,w_n)$,\ $\|w\|=\sqrt{w_1^2+\ldots+w_n^2}$, and let
$$ h_{n}:=\frac1{\|w\|}\,\sum_{j=1}^{n} w_k\,e_j\,. $$
Then
$$ \ag h_{n+1}=\|w\|\,h_{n}+h_{n+2}, \qquad \ag h_{n}=\|w\|\,h_{n+1}. $$
So the Jacobi subspace and Jacobi component of $\gg$ are
\begin{equation}\label{jaccom1}
\cj(\gg)=\lc\{h_j\}_{j\ge n}, \quad J(\gg)=J\bigl(\{0\},\{\|w\|,1,1,\ldots\}\bigr).
\end{equation}

To find the finite-dimensional component, let $\xi=[\xi_{kj}]_{k,j=1}^n$ be a unitary matrix  with the specified last column
\begin{equation}\label{ortmatr}
\xi_{kn}=\frac{w_k}{\|w\|}\,, \quad k=1,\ldots,n.
\end{equation}
We construct an orthonormal basis in $\bc^n$
\begin{equation}\label{finortbas}
f_j:=\sum_{k=1}^n \xi_{kj}e_k(n), \quad k=1,\ldots,n,
\end{equation}
where $\{e_k(n)\}_{k=1}^n$ is the standard basis in $\bc^n$. Put
\begin{equation}\label{star11}
 h_j:=\{f_j,0,0,\ldots\}, \quad j=1,\ldots,n.
\end{equation}
The orthogonality relations $\langle h_k,h_{n}\rangle=0$, \ $1\le k\le n-1$, give
\begin{equation}\label{ortrel}
\ag h_k= \sum_{j=1}^{n} \xi_{kj}\xi_{kn}\cdot h_{n+1}=0.
\end{equation}
Hence the finite-dimensional component $F(\gg)=\mathbb O_{n-1}$ on the subspace
$\cf(\gg)=\lc\{h_j\}_{j=1}^{n-1}$ of the maximal possible dimension $n-1$. So the canonical form is
\begin{equation}\label{canfor1}
\ag\simeq \bo_{n-1}\bigoplus J\bigl(\{0\},\{\|w\|,1,1,\ldots\}\bigr).
\end{equation}

The Jost polynomial is now given by \eqref{jf1}
$$ \|w\|\,u(x)=(1-\|w\|^2)x^2+1. $$
Clearly, $u>0$ for $\|w\|\le1$, and it has zeros inside $(-1,1)$ if and only if $\|w\|>\sqrt2$.
In this case the spectrum is $\s(\gg)=[-2,2]\cup\s_d(\gg)\cup\s_h(\gg)$ with the discrete spectrum being a pair of eigenvalues off $[-2,2]$,
\begin{equation}\label{starweight}
\s_d(S_n(w)+\bp_\infty) =\Bigl\{\pm\Bigl(\sqrt{\|w\|^2-1}+\frac1{\sqrt{\|w\|^2-1}}\Bigr)\Bigr\},
\end{equation}
and the hidden spectrum $\s_h(\gg)=\{0_{n-1}\}$, the zero eigenvalue of multiplicity $n-1$.

For the unweighted star $S_n$ we have
\begin{equation*}
\begin{split}
\s_d(S_n+\bp_\infty) &=\Bigl\{\pm\Bigl(\sqrt{n-1}+\frac1{\sqrt{n-1}}\Bigr)\Bigr\}, \quad \s_h(S_n+\bp_\infty)=0_{n-1},
\quad n\ge3, \\
\s_d(S_2+\bp_\infty) &=\emptyset, \quad \s_h(S_2+\bp_\infty)=0_1.
\end{split}
\end{equation*}
The spectrum $\s(S_n+\bp_\infty)$ was found in \cite{LeNi-umzh}.

Note that the unweighted star graph $S_n$ is a complete bipartite graph, $S_n=K_{1,n}$. For the general complete bipartite graph $K_{p,n+1-p}$
see \cite[Example 5.6]{Gol16}.
\end{example}

Although an explicit form of the matrix $\xi=[\xi_{kj}]_{k,j=1}^n$ in \eqref{ortmatr} is immaterial,
it is worth noting that in the unweighted case  $\xi_{kn}=n^{-1/2}$, $1\le k\le n$, and one can take
\begin{equation}\label{fourmat}
\xi=\Phi_n:=\frac1{\sqrt{n}}\,\bigl[e^{\frac{2\pi ikj}{n}}\bigr]_{k,j=1}^n,
\end{equation}
which is known as the {\it Fourier matrix}. Clearly, there is a number of other options for $\xi$ to be a
real orthogonal matrix (rotation in $\br^n$ with appropriate Euler's angles, orthogonal polynomials etc.).

\begin{example}\label{starmult} ``A multiple star''.

Consider the unweighted star-like graph $S_{n,p}$ with $n$ rays, $n\ge2$, each of which contains
$p+1$ vertices, $p\ge2$. The vertices along each ray (without the root) are numbered as
$$ \{1,n+1,\ldots,(p-1)n+1\},\ \ \{2,n+2,\ldots,(p-1)n+2\},\ \ \ldots \ \ \{n,2n,\ldots,pn\}, $$
and the root is $pn+1$, so $S_{n,1}=S_n$. Let $\gg=S_{n,p}+\bp_\infty$, with the path attached to the root.
As above, we put $h_j:=e_j$, $j=pn+1,\ldots$, and
\begin{equation}
h_{p(n-1)+i}:=\frac1{\sqrt{n}}\,\sum_{q=1}^{n} e_{(i-1)n+q}, \qquad i=1,2,\ldots,p.
\end{equation}
Then $\ag h_{j}=h_{j-1}+h_{j+1}$, $j=pn+2\ldots$, and
\begin{equation*}
\begin{split}
\ag h_{pn+1} &=\sqrt{n}h_{pn}+h_{pn+2}, \ \ \ag h_{pn} =h_{pn-1}+\sqrt{n}\,h_{pn+1}, \\
\ag h_{pn} &=h_{pn-1}+\sqrt{n}\,h_{pn+1}, \\
\ag h_{p(n-1)+i} &=h_{p(n-1)+i-1}+h_{p(n-1)+i+1}, \quad i=2,\ldots,p-1, \\
\ag h_{p(n-1)+1} &=h_{p(n-1)+2},
\end{split}
\end{equation*}
so the Jacobi subspace and Jacobi component of $\gg$ are
\begin{equation}\label{jaccom2}
\begin{split}
\cj(\gg) &=\lc\{h_j\}_{j\ge p(n-1)+1}, \\ J(\gg) &=J(\{0\},\{a_j\}), \quad
a_j=
\left\{
  \begin{array}{rr}
    \sqrt{n}, & j=p; \\
    1, & j\not=p.
  \end{array}
\right.
\end{split}
\end{equation}

To find the finite-dimensional component note that, by the construction,
$h_{p(n-1)+i}\in\lc\{e_{(i-1)n+1},\ldots, e_{in}\}$. As in the above example,
we supplement each $h_{p(n-1)+i}$ to the basis in this subspace by
means of the Fourier matrix \eqref{fourmat}
$$ f_j^{(k)}:=\sum_{q=1}^n \xi_{qj}e_{(k-1)n+q}, \quad f_n^{(k)}:=\sum_{q=1}^n \xi_{qn}e_{(k-1)n+q}
=h_{p(n-1)+k} $$
for $1\le j\le n-1$, $1\le k\le p$. As in \eqref{ortrel} we have
\begin{equation}
\begin{split}\label{adjact}
\ag f_j^{(1)} &=f_j^{(2)}, \quad \ag f_j^{(2)}=f_j^{(1)}+f_j^{(3)}, \ldots, \\
\ag f_j^{(p-1)} &=f_j^{(p-2)}+f_j^{(p)}, \quad \ag f_j^{(p)}=f_j^{(p-1)}.
\end{split}
\end{equation}
Relations \eqref{adjact} mean that the subspace $\ch_j:=\lc\{f_j^{(1)}, \ldots, f_j^{(p)}\}$ is
$A(\gg)$-invariant, and $A(\gg)\vert \ch_j=J_{0,p}$. There are exactly $n-1$ such subspaces
for $j=1,\ldots,n-1$. Finally, we come to the following canonical form for the adjacency matrix
\begin{equation}\label{canfor2}
A(\gg)\simeq F(\gg)\bigoplus J(\gg), \qquad  F(\gg)=\bigoplus_{j=1}^{n-1}J_{0,p}.
\end{equation}

The Jost polynomial is computed in \eqref{rank2}
$$ -\sqrt{n}u(x)=(n-1)x^2\,\frac{x^{2p}-1}{x^2-1}-1=(n-1)(x^{2p}+x^{2p-2}+\ldots+x^2)-1. $$
It is easy to see that $u$ has exactly a pair of symmetric roots $\pm x_0(p,n)$ in $(-1,1)$,
which have the spectral meaning. Hence
$$ \s(\gg)=[-2,2]\cup\s_d(\gg)\cup\s_h(\gg) $$
with
\begin{equation}\label{spec2}
\s_d(\gg)=\Bigl\{\pm\Bigl(x_0(p,n)+\frac1{x_0(p,n)}\Bigr)\Bigr\}\,, \quad \s_h(\gg)=\Bigl\{2\cos\frac{\pi j}{p+1}\Bigr\}_{j=1}^p,
\end{equation}
the hidden spectrum comes from the finite component $F(\gg)$ \eqref{canfor2}, see \eqref{spfinlap}, and each hidden eigenvalue has multiplicity $n-1$.
\end{example}

The problem becomes harder (in the sense of computation) if the original finite star-like graph is
nonsymmetric (the rays are different).

\begin{example}\label{tworays} ``$T_{n,n-1,\infty}$''.

Consider a finite path of order $2n$ with the vertices labeled
$$\{1,2,4,\ldots,2n-2,2n,2n-1,\ldots,5,3\}, $$
with the tail attached to the root $2n$. So there are two finite rays of different length, $n$ and $n-1$, respectively.
The graph $\gg$ thus obtained is known as $T_{n,n-1,\infty}$-graph.

The canonical basis $\{h_j\}_{j\ge1}$ looks as follows: $h_j=e_j$ for $j\ge 2n$,
$$ h_{n+i}=\frac{e_{2i}+e_{2i+1}}{\sqrt2}\,, \quad h_{n-i}=\frac{e_{2i}-e_{2i+1}}{\sqrt2}\,, \qquad i=1,\ldots,n-1, $$
and $h_n=e_1$.

The adjacency operator $A(\gg)$ acts on the basis vectors in a simple way:
\begin{equation*}
\begin{split}
\ag h_{2n+j} &=h_{2n+j+1}+h_{2n+j-1}, \quad j=1,2,\ldots, \\
\ag h_{2n} &=h_{2n+1}+\sqrt2\,h_{2n-1}, \qquad \ag h_{2n-1} =\sqrt2\,h_{2n}+h_{2n-2}, \\
\ag h_{2n-k} &=h_{2n-k+1}+h_{2n-k-1}, \quad k=2,\ldots,n-2.
\end{split}
\end{equation*}

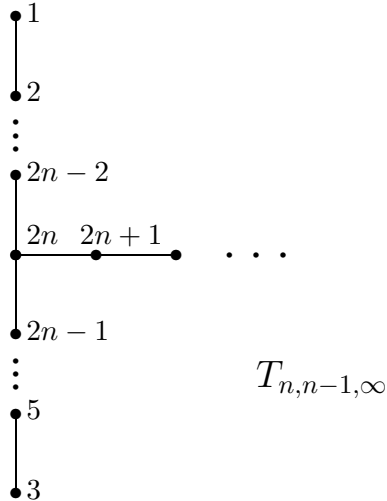
\begin{picture}(300, 200)

\multiput(60,100) (30,0) {3} {\circle* {4}}
\multiput(140,100) (10,0) {3} {\circle* {2}}
\multiput(60,190) (0,-30) {2} {\circle* {4}}
\multiput(60,150) (0,-5) {3} {\circle* {2}}
\multiput(60,130) (0,-30) {3} {\circle* {4}}
\multiput(60,60) (0,-5) {3} {\circle* {2}}
\multiput(60,40) (0,-30) {2} {\circle* {4}}
\multiput(62,100) (30,0) {2} {\line(1, 0) {26}}
\multiput(60,128) (0,-30) {2} {\line(0,-1) {26}}
\multiput(60,188) (0,-30) {1} {\line(0,-1) {26}}
\multiput(60,38) (0,-30) {1} {\line(0,-1) {26}}
\put(64,104) {$2n$}
\put(84,104) {$2n+1$}
\put(64,08) {$3$}
\put(64,38) {$5$}
\put(64,68) {$2n-1$}
\put(64,188) {$1$}
\put(64,158) {$2$}
\put(64,128) {$2n-2$}
\put(150,50) {\bf  \Large{$T_{n,n-1,\infty}$}}
\end{picture}

Furthermore,
\begin{equation*}
\begin{split}
\ag h_{n+1} &=\frac{e_1+e_4+e_5}{\sqrt2}=h_{n+2}+\frac1{\sqrt2}\,h_n \\
\ag h_{n} &=e_2=\frac{h_{n+1}+h_{n-1}}{\sqrt2}, \qquad \ag h_{n-j} =h_{n-j+1}+h_{n-j-1}, \quad j=2,\ldots,n-2, \\
\ag h_{n-1} &=\frac{e_1+e_4-e_5}{\sqrt2}=\frac1{\sqrt2}\,h_n+h_{n-2}, \quad \ag h_1=h_2.
\end{split}
\end{equation*}
So, the finite-dimensional component is missing, and $\ag\simeq J(\{0\}, \{a_j\})$,
$$ a_j=1, \ \ j\not=n-1,n,2n-1, \quad a_{n-1}=a_n=\frac1{\sqrt2}\,, \quad a_{2n-1}=\sqrt2. $$

The Jost polynomial can be found from the definition \eqref{3term}. For instance, for $n=3$ we have
$$ \frac1{\sqrt2}\,u(x)=-ax^{10}-2ax^8-\Bigl(\frac52 a-1\Bigr)x^6-\Bigl(2a+\frac32\Bigr)x^4-ax^2+1, \quad a=\sqrt2-1. $$
Then $u$ has a pair of symmetric roots $\pm x_1$ in $(-1,1)$, and so
\begin{equation}\label{oneeig}
\s(T_{3,2,\infty})=[-2,2]\cup\s_d(T_{3,2,\infty}), \qquad \s_d(T_{3,2,\infty})=\pm\Bigl(x_1+\frac1{x_1}\Bigr).
\end{equation}
\end{example}

\begin{example}\label{sword} ``$T_{1,1,2,\infty}$''.

\begin{picture}(300, 130)
\multiput(60,80) (30,0) {5} {\circle* {4}}
\multiput(200,80) (10,0) {3} {\circle* {2}}
\multiput(120,110) (0,-30) {3} {\circle* {4}}
\multiput(62,80) (30,0) {4} {\line(1, 0) {26}}
\multiput(120,108) (0,-30) {2} {\line(0,-1) {26}}
\put(60,84) {$1$}
\put(94,84) {$2$}
\put(122,83) {$5$}
\put(152,83) {$6$}
\put(182,83) {$7$}
\put(124,110) {$3$}
\put(124,50) {$4$}
\put(210,30) {\bf  \Large{$T_{1,1,2,\infty}$}}
\end{picture}

This graph can be viewed as a coupling $\gg=T(1,2,2)+\bp_\infty$. To construct the canonical basis we put $h_j=e_j$ for $j\ge5$, and
$$ h_4=\frac{e_2+e_3+e_4}{\sqrt3}\,, \ \ h_3=e_1, \ \ h_2=\sqrt{\frac23}\,e_2-\frac{e_3+e_4}{\sqrt6}\,. $$
Then $\cj(\gg)=\lc\{h_j\}_{j\ge2}$, and the adjacency operator acts as
\begin{equation*}
\begin{split}
\ag h_j &=h_{j+1}+h_{j-1}, \quad j=6,7,\ldots, \\
\ag h_5 &=h_{6}+\sqrt3\,h_{4}, \quad \ag h_4 =\sqrt3\,h_{5}+\frac1{\sqrt3}\,h_{3}, \\
\ag h_3 &= \frac1{\sqrt3}\,h_{4}+\sqrt{\frac23}\,h_{2}, \quad \ag h_2=\sqrt{\frac23}\,h_{3}.
\end{split}
\end{equation*}

If we add a vector $h_1=\frac{e_3-e_4}{\sqrt2}$, then $\ag\,h_1=0$, so the finite-dimensional component has dimension $1$ and $F(\gg)=0$. Finally,
\begin{equation}
\ag\simeq \bo_{1}\bigoplus J\Bigl(\{0\},\Bigl\{\sqrt{\frac23}, \frac1{\sqrt3}, \sqrt3, 1,1,\ldots \Bigr\}\Bigr).
\end{equation}

The Jost polynomial for the Jacobi component
$$ \sqrt{\frac23}\,u(x)=-2x^6-2x^4-x^2+1 $$
has a pair of symmetric roots $\pm x_2$ in $(-1,1)$, so
$$ \s(T_{1,1,2,\infty})=[-2,2]\cup\s_d(T_{1,1,2,\infty})\cup\s_h(T_{1,1,2,\infty}) $$
with
\begin{equation*}
\s_d(T_{1,1,2,\infty})=\Bigl\{\pm\Bigl(x_2+\frac1{x_2}\Bigr)\Bigr\}, \quad \s_h(T_{1,1,2,\infty})=0_1.
\end{equation*}

The similar example of the graph $T_{1,2,2,\infty}$ is studied in \cite[Example 4.4]{Gol16}.
\end{example}

\begin{example}\label{tree}  ``A tree with tail''.

Consider the coupling of a simple tree graph of order $8$ with the infinite path. The canonical basis looks as follows: $h_j=e_j$ for $j\ge8$,
\begin{equation*}
\begin{split}
h_7 &=\frac{e_6+e_7}{\sqrt2}\,, \quad h_6=\frac1{\sqrt5}\,\sum_{i=1}^5 e_i, \quad h_5 =\frac{e_6-e_7}{\sqrt2}\,, \\
h_4 &=\sqrt{\frac2{15}}\,(e_1+e_2+e_3)-\sqrt{\frac3{10}}\,(e_4+e_5).
\end{split}
\end{equation*}

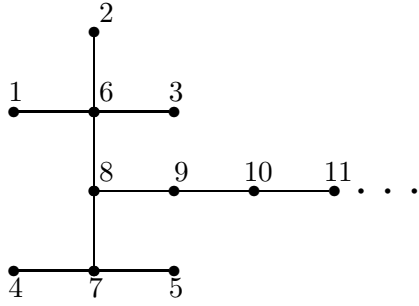
\begin{picture}(300, 140)
\multiput(90,60) (30,0) {4} {\circle* {4}}
\multiput(190,60) (10,0) {3} {\circle* {2}}
\multiput(90,120) (0, -30) {4} {\circle* {4}}
\multiput(60,90) (30,0) {3} {\circle* {4}}
\multiput(60,30) (30,0) {3} {\circle* {4}}
\multiput(62,30) (30,0) {2} {\line(1, 0) {26}}
\multiput(92,60) (30,0) {3} {\line(1, 0) {26}}
\multiput(62,90) (30,0) {2} {\line(1, 0) {26}}
\multiput(90,118) (0,-30) {3} {\line(0,-1) {26}}
\put(58,94) {$1$}
\put(58,20) {$4$}
\put(88,20) {$7$}
\put(118,20) {$5$}
\put(92,94) {$6$}
\put(118,94) {$3$}
\put(92,124) {$2$}
\put(92,64) {$8$}
\put(120,64) {$9$}
\put(146,64) {$10$}
\put(176,64) {$11$}
\end{picture}

The adjacency operator acts on the vectors as
\begin{equation*}
\begin{split}
\ag h_8 &=h_9+\sqrt2\,h_7, \quad \ag h_7 =\sqrt2\,h_8+\sqrt{\frac52}\,h_6, \\
\ag h_6 &=\sqrt{\frac52}\,h_7+\frac1{\sqrt{10}}\,h_5, \quad \ag h_5=\frac1{\sqrt{10}}\,h_6+2\sqrt{\frac35}\,h_4, \\
\ag h_4 &=2\sqrt{\frac35}\,h_5.
\end{split}
\end{equation*}
Hence, the Jacobi subspace is $\cj(T)=\lc\{h_j\}_{j\ge4}$, and the Jacobi component is
$$ J(\gg)=J\Bigl(\{0\},\Bigl\{2\sqrt{\frac35}, \frac1{\sqrt{10}}, \sqrt{\frac52}, \sqrt2, 1, 1,\ldots\Bigr\}\Bigr). $$

To find the finite-dimensional component, we supplement the sequence $\{h_j\}_{j\ge4}$ with
$$ h_3=\frac{e_4-e_5}{\sqrt2}\,, \quad h_2=\frac{e_1-e_3}{\sqrt2}\,, \quad h_1=\frac{e_1-2e_2+e_3}{\sqrt6} $$
and note that $\ag h_3=\ag h_2=\ag h_1=0$, so $F(\gg)=\bo_3$.

The Jost polynomial is
$$ \sqrt{\frac65}\,u(x)=(x^2+1)(-x^6+x^4+2x^2-1), $$
again with a pair of symmetric roots $\pm x_3$ in $(-1,1)$. So, the spectrum is
\begin{equation*}
\s(\gg)=[-2,2]\cup\Bigl\{\pm\Bigl(x_3+\frac1{x_3}\Bigr)\Bigr\}\cup\{0_3\}.
\end{equation*}
\end{example}

\subsection{Graphs with cycles and tails}

\begin{example}\label{3-flower}  ``A flower with equal petals of degree $3$''.

Consider $n\ge2$ copies of the cycle $\bc_3$ of order $3$, glued together at one common vertex (root) $2n+1$, so that $j$-th copy is $\{2j-1,2j,2n+1\}$, $j=1,2,\ldots,n$.
We denote this graph by $(\bc_3)^n$.
Let $\gg=(\bc_3)^n+\bp_\infty$ be the coupling of this graph with the infinite path $\{2n+1, 2n+2,\ldots\}$ attached to the root.

To construct the canonical basis, put
$$ h_k=e_k, \quad k\ge 2n+1, \quad h_{2n}=\frac1{\sqrt{2n}}\,\sum_{i=1}^{2n} e_i. $$
Then $\ag h_{2n+1}=h_{2n+2}+\sqrt{2n}\,h_{2n}$, so the Jacobi component is
$$ \cj(\gg)=\lc\{h_j\}_{j\ge 2n}, \quad J(\gg)=J(\{1,0,0,\ldots\}, \{\sqrt{2n},1,1,\ldots\}). $$

Next, let $\xi=[\xi_{ij}]$ be a unitary $n\times n$-matrix with $\xi_{i1}=n^{-1/2}$ (the constant first column). We define a matrix $\eta=[\eta_{ij}]$ of order $2n$
by the following recipe
\begin{equation*}
\begin{split}
\eta_{2k-1,j} &=\eta_{2k,j}=\xi_{kj}/\sqrt2, \quad k=1,2,\ldots,n, \quad j=1,2,\ldots,n; \\
\eta_{2k-1,j} &=-\eta_{2k,j}=\xi_{k,j-n}/\sqrt2, \quad k=1,2,\ldots,n, \quad j=n+1,n+2,\ldots,2n.
\end{split}
\end{equation*}
The matrix $\eta$ is unitary, with the constant first column $\eta_{i1}=(2n)^{-1/2}$. Put
$$ h_{2n-j}=\sum_{i=1}^{2n} \eta_{i,j+1}e_j, \qquad j=0,1,\ldots,2n-1. $$
The adjacency operator $\ag$ act on the basis vectors as
\begin{equation*}
\begin{split}
\ag h_{2n-j} &=\sum_{i=1}^{2n} \eta_{i,j+1}\,\ag e_i=\sum_{k=1}^{n} \eta_{2k,j+1}\,\ag e_{2k}+\sum_{k=1}^{n} \eta_{2k-1,j+1}\,\ag e_{2k-1} \\
{} &=\sum_{i=1}^{2n} \eta_{i,j+1}\,e_{2n+1}+\sum_{k=1}^{n} \eta_{2k,j+1}\,e_{2k-1}+\sum_{k=1}^{n} \eta_{2k-1,j+1}\,e_{2k}.
\end{split}
\end{equation*}
Finally,
$$ \ag h_{2n-j}=\left\{
                    \begin{array}{ll}
                      h_{2n-j}, & \hbox{$j=1,2,\ldots,n-1$;} \\
                      -h_{2n-j}, & \hbox{$j=n,n+1,\ldots,2n-1$.}
                    \end{array}
                  \right.
$$
The finite-dimensional component is $F(\gg)={\rm diag}(-1,-1,\ldots,-1,1,1,\ldots,1)$ with $n$ numbers $-1$ and $n-1$ numbers $1$. So,
$$ \ag\simeq {\rm diag}(-1,-1,\ldots,-1,1,1,\ldots,1)\bigoplus J(\{1,0,0,\ldots\}, \{\sqrt{2n},1,1,\ldots\}). $$

The Jost polynomial, computed in \eqref{jf1},
$$ \sqrt{2n}\,u(x)=(1-2n)x^2-x+1, $$
has two roots $-1<x_5<0<x_4<1$,
$$ x_{4,5}=\frac{-1\pm\sqrt{1+4(2n-1)}}{2(2n-1)}\,. $$
Hence,
\begin{equation*}
\s(\gg)=[-2,2]\cup\Bigl\{x_4+\frac1{x_4}, \ x_5+\frac1{x_5}\Bigr\}\cup\{(-1)_{n}, \ 1_{n-1}\}.
\end{equation*}

For the spectrum of a propeller with {\it two} equal blades of an arbitrary order see \cite[Example 5.2]{Gol16}. In Example \ref{flower} below
we study the general flower graphs with the infinite tail.
\end{example}

\begin{example}\label{complete} ``The complete graph with tail''.

Let $K_n$ be a complete graph of order $n\ge3$, $\gg=K_n+\bp_\infty$, the ray $\{n,n+1,\ldots\}$ is attached to the vertex $n$. Put
$$ h_j=e_j, \quad j=n,n+1\ldots, \quad h_{n-1}=\frac1{\sqrt{n-1}}\,\sum_{k=1}^{n-1} e_k. $$
Since
$$ \ag e_k=\sum_{j\not=k} e_j=S-e_k, \qquad S:=\sum_{j=1}^n e_j, \quad k=1,2,\ldots,n-1, $$
we see that
\begin{equation*}
\begin{split}
\ag h_{n-1} &=\frac1{\sqrt{n-1}}\,\sum_{k=1}^{n-1} (S-e_k)=\sqrt{n-1}\,S-\whe_{n-1} \\
{} &=\sqrt{n-1}\bigl(\sqrt{n-1}\,h_{n-1}+h_n\bigr)-h_{n-1}=\sqrt{n-1}\,h_n+(n-2)\,h_{n-1}.
\end{split}
\end{equation*}
Hence, the Jacobi component is
$$ \cj(\gg)=\lc\{h_j\}_{j\ge n-1}, \quad J(\gg)=J(\{n-2,0,0,\ldots\}, \{\sqrt{n-1},1,1,\ldots\}). $$

Next, put
$$ h_k=\sum_{j=1}^{n-2}\xi_{kj}e_j, \quad \xi_{k,n-2}=\frac1{\sqrt{n-2}}, \quad k=1,2,\ldots,n-2, $$
$[\xi_{ij}]$ is a unitary matrix of order $n-2$. Then,
$$ \ag h_k=\sum_{j=1}^{n-2}\xi_{kj}(S-e_j)=S\,\sum_{j=1}^{n-2}\xi_{kj}-h_k=-h_k, \quad k=1,\ldots,n-2. $$
So, the finite dimensional component is
$$ \cf(\gg)=\lc\{h_j\}_{j=1}^{n-2}, \quad F(\gg)=-I_{n-2}. $$

The Jost polynomial, computed in \eqref{jf1}
$$ \sqrt{n-1}\,u(x)=-(n-2)x^2-(n-2)x+1, $$
has one root $x_6$ in $(-1,1)$,
$$ x_{6}=\frac12\Bigl(\sqrt{\frac{n+2}{n-2}}-1\Bigr)\,. $$
Hence,
\begin{equation*}
\s(\gg)=[-2,2]\cup\Bigl\{x_6+\frac1{x_6}\Bigr\}\cup\{(-1)_{n-2}\}.
\end{equation*}
\end{example}

\begin{example}\label{flag}  ``A flag with a tail''.

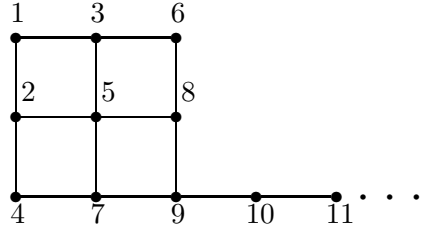
\begin{picture}(300, 120)
\multiput(60,30) (30,0) {5} {\circle* {4}}
\multiput(190,30) (10,0) {3} {\circle* {2}}
\multiput(60,60) (30,0) {3} {\circle* {4}}
\multiput(60,90) (30,0) {3} {\circle* {4}}
\multiput(62,30) (30,0) {4} {\line(1, 0) {26}}
\multiput(62,60) (30,0) {2} {\line(1, 0) {26}}
\multiput(62,90) (30,0) {2} {\line(1, 0) {26}}
\multiput(60,88) (0,-30) {2} {\line(0,-1) {26}}
\multiput(90,88) (0,-30) {2} {\line(0,-1) {26}}
\multiput(120,88) (0,-30) {2} {\line(0,-1) {26}}
\put(58,20) {$4$}  \put(88,20) {$7$}
\put(118,20) {$9$}  \put(146,20) {$10$} \put(176,20) {$11$}
\put(58,96) {$1$}  \put(88,96) {$3$} \put(118,96) {$6$}
\put(62,66) {$2$}  \put(92,66) {$5$} \put(122,66) {$8$}
\end{picture}

The canonical basis looks as follows: $h_j=e_j$ for $j=9,10,\ldots$,
$$ h_8=\frac{e_8+e_7}{\sqrt2}\,, \ \ h_7=\frac{e_6+2e_5+e_4}{\sqrt6}\,, \ \ h_6=\frac{e_3+e_2}{\sqrt2}\,, \ \ h_5=e_1. $$
The Jacobi component is
$$ \cj(\gg)=\lc\{h_j\}_{j\ge5}, \quad J(\gg)=J(\{0\}, \{\sqrt2,\sqrt3,\sqrt3,\sqrt2,1,1,\ldots\}). $$

To complete the basis, we put
$$ h_4=\frac{e_8-e_7}{\sqrt2}\,, \ \ h_3=\frac{e_6-e_4}{\sqrt2}\,, \ \ h_2=\frac{e_3-e_2}{\sqrt2}\,, \ \ h_1=\frac{e_6-e_5+e_4}{\sqrt3}. $$
It is easy to see that
$$ \ag h_4=h_3, \ \ \ag h_3=h_4-h_2, \ \ \ag h_2=h_3, \ \ \ag h_1=0. $$
So, the finite-dimensional component is
$$ F(\gg)=\begin{bmatrix}
0 & 0 & 0 & 0 \\
0 & 0 & -1 & 0  \\
0 & 1 & 0 & 1 \\
0 & 0 & 1 & 0
\end{bmatrix}, \quad \s(F(\gg))=0_4.
$$

The Jost polynomial
$$ 6u(x)=-x^8-6x^2+1 $$
has a pair of symmetric roots $\pm x_7$ in $(-1,1)$, so
\begin{equation*}
\s(\gg)=[-2,2]\cup\Bigl\{\pm\Bigl(x_7+\frac1{x_7}\Bigr)\Bigr\}\cup\{0_4\}.
\end{equation*}
\end{example}

\begin{example}\label{cycle2} ``A cycle with two tails''.

Let $\bc_m$ be a cycle of order $m$. The spectral analysis of the coupling $\bc_m+\bp_\infty$ was carried out in \cite[Proposition 1.4]{Gol16}.
We consider here the graph $\gg=\bc_{2n+1}+\bp_\infty+\bp_\infty$, $n\ge2$, with two tails,
$$ \{2n, 2n+2\ldots\}, \quad \{2n+1,2n+3,\ldots\}, $$
attached to the adjacent vertices $2n$ and $2n+1$, respectively. Although $\gg$ is not a coupling in the sense of Definition \ref{coupl},
the methods works in this situation as well.

Define a system of vectors $h_0=e_n$,
\begin{equation*}
\begin{split}
h_k^{\pm} &=\frac{e_{n+k}\pm e_{n-k}}{\sqrt2}\,, \qquad k=1,2,\ldots,n-1, \\
h_{n+i}^{\pm} &=\frac{e_{2n+2i}\pm e_{2n+2i+1}}{\sqrt2}\,, \qquad i=0,1,\ldots
\end{split}
\end{equation*}
The system $\{h_0,h_j^\pm\}_{j\ge1}$ is the canonical orthonormal basis. The adjacency operator $\ag$ acts as
\begin{equation*}
\begin{split}
\ag h_0 &=\sqrt2\,h_1^+, \quad \ag h_1^+=h_2^+ +\sqrt2\,h_0, \\
\ag h_k^+ &=h_{k+1}^+ +h_{k-1}^+, \quad k=2,\ldots,n-1, \\
\ag h_n^+ &=h_{n+1}^+ +h_{n}^+ +h_{n-1}^+, \\
\ag h_k^+ &=h_{k+1}^+ +h_{k-1}^+, \quad k\ge n+1.
\end{split}
\end{equation*}
So, the first Jacobi component is
$$ \cj^+(\gg)=\lc\{h_0, \{h_j^+\}\}_{j\ge1}, \quad J^+(\gg)=J(\{b_j^+\}, \{a_j^+\}) $$
with
$$ b_j^+=\left\{
           \begin{array}{rr}
             0, & \hbox{$j\not=n+1$;} \\
             1, & \hbox{$j=n+1$,}
           \end{array}
         \right.  \qquad
   a_j^+=\left\{
           \begin{array}{rr}
             1, & \hbox{$j\not=1$;} \\
             \sqrt2, & \hbox{$j=1$.}
           \end{array}
         \right.
$$

Next,
\begin{equation*}
\begin{split}
\ag h_1^- &=h_2^-, \quad \ag h_k^- =h_{k+1}^- +h_{k-1}^-, \quad k=2,\ldots,n-1, \\
\ag h_n^- &=h_{n+1}^- -h_{n}^- +h_{n-1}^-, \\
\ag h_k^- &=h_{k+1}^- +h_{k-1}^-, \quad k\ge n+1.
\end{split}
\end{equation*}
So, the second Jacobi component is
$$ \cj^-(\gg)=\lc\{h_j^-\}_{j\ge1}, \quad J^-(\gg)=J(\{b_j^-\}, \{a_j^-\}) $$
with
$$ b_j^-=\left\{
           \begin{array}{rr}
             0, & \hbox{$j\not=n$;} \\
             -1, & \hbox{$j=n$,}
           \end{array}
         \right.  \qquad
   a_j^-\equiv 1, \qquad j=1,2,\ldots
$$
The canonical form of the adjacency operator is the orthogonal sum of two Jacobi operators
$$ \ag\simeq J^+(\gg)\bigoplus J^-(\gg). $$

The Jost polynomial for $J^+(\gg)$, given in \eqref{jf4} with $q=n+1$,
$$ \sqrt2 u_+(x)=-x^{2n+1}-x^2-x+1, $$
has one root $x_8>0$ in $(-1,1)$. Similarly, the Jost polynomial for $J^-(\gg)$, given in \eqref{rank1},
$$ -u(x)=1+x\frac{x^{2n}-1}{x^2-1}, $$
has one root $x_9<0$ in $(-1,1)$. Hence,
\begin{equation*}
\s(\gg)=[-2,2]\cup\s_d(\gg), \quad \s_d(\gg)=\Bigl\{x_8+\frac1{x_8}, x_9+\frac1{x_9}\Bigr\}.
\end{equation*}

The case of the cycle $\bc_{2n+1}$ with two tails attached to the {\it same} vertex is studied in \cite[Example 5.1]{Gol16}.
\end{example}

\section{Ladders and chains of cycles}
\label{s4}

\subsection{Canonical form and periodic structure}

For certain infinite graphs (not necessarily finite graphs with tails) the canonical basis arises in a natural way. The canonical form, distinct from
\eqref{canform}, enables one to describe explicitly the spectrum of the graph in question.

\begin{example}\label{compladd} ``A complete ladder''.

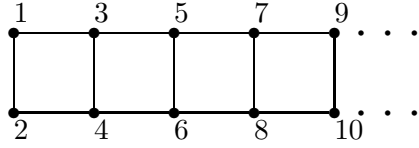
\begin{picture}(300, 90)
\multiput(60,30) (30,0) {5} {\circle* {4}}
\multiput(60,60) (30,0) {5} {\circle* {4}}
\multiput(62,30) (30,0) {4} {\line(1, 0) {26}}
\multiput(62,60) (30,0) {4} {\line(1, 0) {26}}
\multiput(190,30) (10,0) {3} {\circle* {2}}
\multiput(190,60) (10,0) {3} {\circle* {2}}
\multiput(60,58) (0,-30) {1} {\line(0,-1) {26}}
\multiput(90,58) (0,-30) {1} {\line(0,-1) {26}}
\multiput(120,58) (0,-30) {1} {\line(0,-1) {26}}
\multiput(150,58) (0,-30) {1} {\line(0,-1) {26}}
\multiput(180,58) (0,-30) {1} {\line(0,-1) {26}}
\put(60,20) {$2$}  \put(90,20) {$4$}  \put(120,20) {$6$}  \put(150,20) {$8$}  \put(180,20) {$10$}
\put(60,64) {$1$}  \put(90,64) {$3$}  \put(120,64) {$5$}  \put(150,64) {$7$}  \put(180,64) {$9$}
\end{picture}

We define the canonical basis $\{h^+_n, \ h_n^-\}_{n\ge1}$ by the relations
$$ h_n^\pm:=\frac{e_{2n-1}\pm e_{2n}}{\sqrt2}, \qquad n=1,2,\ldots. $$
The subspaces $\cj^\pm=\lc\{h_n^\pm\}_{n\ge1}$ are invariant for $\ag$, $\ell^2=\cj^+\oplus\cj^-$, and the restrictions
on these subspaces $J^\pm(\gg)=\pm I+J_0$. So, the canonical form of $\gg$ is
\begin{equation*}
\ag\simeq\begin{bmatrix}
I+J_0 &  \\
  & -I+J_0
\end{bmatrix}.
\end{equation*}
The spectrum is $\s(\gg)=[-3,3]$, with the interval $[-1,1]$ of multiplicity $2$.

Similarly, for the ladder below

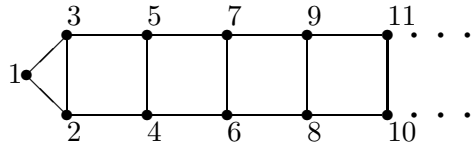
\begin{picture}(300, 90)
\multiput(60,30) (30,0) {5} {\circle* {4}}
\multiput(60,60) (30,0) {5} {\circle* {4}}
\multiput(62,30) (30,0) {4} {\line(1, 0) {26}}
\multiput(62,60) (30,0) {4} {\line(1, 0) {26}}
\multiput(190,30) (10,0) {3} {\circle* {2}}
\multiput(190,60) (10,0) {3} {\circle* {2}}
\multiput(60,58) (0,-30) {1} {\line(0,-1) {26}}
\multiput(90,58) (0,-30) {1} {\line(0,-1) {26}}
\multiput(120,58) (0,-30) {1} {\line(0,-1) {26}}
\multiput(150,58) (0,-30) {1} {\line(0,-1) {26}}
\multiput(180,58) (0,-30) {1} {\line(0,-1) {26}}
\put(60,20) {$2$}  \put(90,20) {$4$}  \put(120,20) {$6$}  \put(150,20) {$8$}  \put(180,20) {$10$}
\put(60,64) {$3$}  \put(90,64) {$5$}  \put(120,64) {$7$}  \put(150,64) {$9$}  \put(180,64) {$11$}
\put(45,45) {\line(1,1){15}}  \put(45,45) {\line(1,-1){15}}
\put(45,45) {\circle*{4}}
\put(38,42) {$1$}
\end{picture}

we define
$$ h_n^\pm:=\frac{e_{2n}\pm e_{2n+1}}{\sqrt2}, \quad n=1,2,\ldots, \quad h_0^+=e_1. $$
The adjacency operator $\ag$ acts as
\begin{equation*}
\begin{split}
\ag h_0^+ &=\sqrt2\,h_1^+, \quad \ag h_1^+=\sqrt2\,h_0^+ +h_1^+ +h_2^+, \\
\ag h_n^+ &=h_{n-1}^+ + h_n^+ + h_{n+1}^+, \quad n=2,3,\ldots,
\end{split}
\end{equation*}
\begin{equation*}
\begin{split}
\ag h_1^- &=-h_1^- +h_2^-, \\
\ag h_n^- &=h_{n-1}^- -h_n^- +h_{n+1}^-, \quad n=2,3,\ldots.
\end{split}
\end{equation*}

Again, the subspaces $\cj^+=\lc\{h_n^+\}_{n\ge0}$ and $\cj^-=\lc\{h_n^-\}_{n\ge1}$ are invariant for $\ag$, $\ell^2=\cj^+\oplus\cj^-$,
and the canonical form of $\gg$ is
\begin{equation}\label{canlad}
\ag\simeq\begin{bmatrix}
I+J^+ &  \\
  & -I+J^-
\end{bmatrix},
\end{equation}
where
$$ J^-=J_0, \qquad J^+=J(\{-1,0,0,\ldots\}, \{\sqrt2,1,1,\ldots\}). $$

The Jost polynomial for $J^+$ is computed in \eqref{jf1}
$$ \sqrt2\,u(x)=-x^2+x+1, $$
with the roots
$$ x_{10,11}=\frac{1\pm\sqrt5}2\,, $$
so $x_{11}\in(-1,1)$. Finally, $\s(I+J^+)=[-1,3]\cup \{1-\sqrt5\}$, and
$$ \s(\gg)=[-3,3]\cup\s_h(\gg), \quad \s_h(\gg)=\{1-\sqrt5\}, $$
with the eigenvalue $1-\sqrt5\in[-2,-1]$ lying on the simple spectrum.
\end{example}

\begin{example}\label{lantern} ``A lantern on a ladder''.

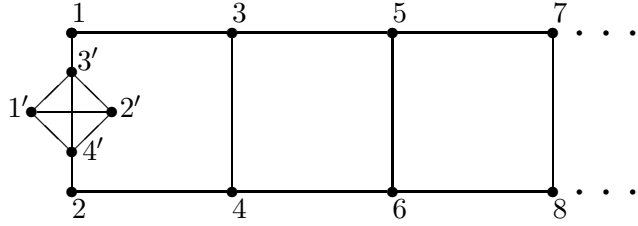
\begin{picture}(300, 120)
\multiput(60,30) (60,0) {4} {\circle* {4}}
\multiput(60,90) (60,0) {4} {\circle* {4}}
\multiput(250,30) (10,0) {3} {\circle* {2}}
\multiput(250,90) (10,0) {3} {\circle* {2}}
\multiput(60,90) (0,-15) {1} {\circle* {4}}
\multiput(60,75) (0,-30) {1} {\circle* {4}}
\multiput(60,45) (0,-15) {1} {\circle* {4}}
\multiput(120,88) (0,-60) {1} {\line(0,-1) {56}}
\multiput(180,88) (0,-60) {1} {\line(0,-1) {56}}
\multiput(240,88) (0,-60) {1} {\line(0,-1) {56}}
\multiput(62,30) (60,0) {3} {\line(1,0) {56}}
\multiput(62,90) (60,0) {3} {\line(1,0) {56}}
\multiput(60,88) (0,-15) {1} {\line(0,-1) {11}}
\multiput(60,73) (0,-30) {1} {\line(0,-1) {26}}
\multiput(60,43) (0,-15) {1} {\line(0,-1) {11}}
\put(45,60) {\circle*{4}}  \put(75,60) {\circle*{4}}
\multiput(47,60) (30,0) {1} {\line(1,0) {26}}
\put(45,60) {\line(1,1){15}} \put(60,45) {\line(1,1){15}}
\put(45,60) {\line(1,-1){15}} \put(60,75) {\line(1,-1){15}}
\put(60,94) {$1$}  \put(120,94) {$3$}  \put(180,94) {$5$}  \put(240,94) {$7$}
\put(60,20) {$2$}  \put(120,20) {$4$}  \put(180,20) {$6$}  \put(240,20) {$8$}
\put(36,58) {$1'$}  \put(78,58) {$2'$}  \put(62,76) {$3'$}  \put(64,42) {$4'$}
\end{picture}

The standard basis for the Hilbert space $\ell^2(\gg)$ of this graph is enumerated as $\{e_1', e_2', e_3', e_4'; e_k, k\ge1\}$.
The canonical basis is
$$ \{g_1^+, g_2^+, g_1^-, g_2^-; h_j^+, h_j^-, j\ge1\} $$
with
$$ g_j^\pm=\frac{e_{2j-1}'\pm e_{2j}'}{\sqrt2}\,, \quad j=1,2, \qquad h_n^\pm=\frac{e_{2n-1}\pm e_{2n}}{\sqrt2}\,, \quad n=1,2,\ldots. $$
The adjacency operator $\ag$ acts as
\begin{equation*}
\begin{split}
\ag g_1^+ &=g_1^+ +2g_2^+, \quad \ag g_2^+=2g_1^+ +g_2^+ +h_1^+, \\
\ag h_1^+ &=g_2^+ +h_2^+, \quad \ag h_k^+=h_{k-1}^+ +h_k^+ +h_{k+1}^+, \quad k=2,3,\ldots
\end{split}
\end{equation*}
and
\begin{equation*}
\begin{split}
\ag g_1^- &=-g_1^-, \quad \ag g_2^-=-g_2^- +h_1^-, \\
\ag h_1^- &=g_2^- +h_2^-, \quad \ag h_k^-=h_{k-1}^- -h_k^- +h_{k+1}^-, \quad k=2,3,\ldots.
\end{split}
\end{equation*}
So, we have the following decomposition of the space $\ell^2(\gg)$ in orthogonal sum of three $\ag$-invariant subspaces
$$ \ell^2(\gg)=\cj'\oplus\cj^+\oplus\cj^-, \qquad \ag=J'\oplus J^+\oplus J^-, $$
where
$$ \cj'=\lc\{g_1^-\}, \quad \cj^+=\lc\{g_1^+, g_2^+; h_k^+, k\ge1\}, \quad \cj^-=\lc\{g_2^-; h_k^-, k\ge1\}. $$
The restrictions of $\ag$ on these subspaces look $J'=-I_1$,
$$ J^+=I+J(\{0,0,-1,0,0,\ldots\}, \{2,1,1,\ldots\}), \quad J^-=-I+J(\{0,1,0,0,\ldots\}, \{1\}). $$

Both the matrices above are of finite rank. The Jost polynomial for the second one is given in \eqref{rank1} with $q=2$
$$ u(x)=1-x\,\frac{x^4-1}{x^2-1}=-x^3-x^2-x+1. $$
It has the only root $x_{12}>0$ in $(-1,1)$, so
\begin{equation}\label{lant1}
\s(J^-)=[-3,1]\cup\,\Bigl\{x_{12}+\frac1{x_{12}}-1\Bigr\}.
\end{equation}

Relation \eqref{jf4} with $q=3$ provides the Jost polynomial for the first matrix
$$ 2u(x)=\frac{x^7-3x^5+3x^3-x}{x^2-1}-3x^2+1=x^5-2x^3-3x^2+x+1. $$
By using the Descarte's rule, one can see that $u$ has two roots $x_{13}<0<x_{14}$ in $(-1,1)$, so
\begin{equation}\label{lant2}
\s(J^+)=[-1,3]\cup\,\Bigl\{x_{13}+\frac1{x_{13}}+1, x_{14}+\frac1{x_{14}}+1\Bigr\}.
\end{equation}

The union of \eqref{lant1}, \eqref{lant2}, and an eigenvalue $(-1)_1$, constitutes $\s(\gg)$.
\end{example}

\begin{example}\label{simlad} ``Simon's ladder'' \cite{Si96}.

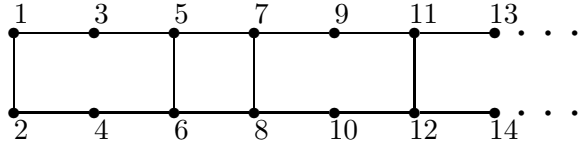
\begin{picture}(300, 90)
\multiput(60,30) (30,0) {7} {\circle* {4}}
\multiput(60,60) (30,0) {7} {\circle* {4}}
\multiput(62,30) (30,0) {6} {\line(1, 0) {26}}
\multiput(62,60) (30,0) {6} {\line(1, 0) {26}}
\multiput(250,30) (10,0) {3} {\circle* {2}}
\multiput(250,60) (10,0) {3} {\circle* {2}}
\multiput(60,58) (0,-30) {1} {\line(0,-1) {26}}
\multiput(120,58) (0,-30) {1} {\line(0,-1) {26}}
\multiput(150,58) (0,-30) {1} {\line(0,-1) {26}}
\multiput(210,58) (0,-30) {1} {\line(0,-1) {26}}
\put(60,20) {$2$}  \put(90,20) {$4$}  \put(120,20) {$6$}  \put(150,20) {$8$}  \put(178,20) {$10$}  \put(208,20) {$12$}  \put(238,20) {$14$}
\put(60,64) {$1$}  \put(90,64) {$3$}  \put(120,64) {$5$}  \put(150,64) {$7$}  \put(180,64) {$9$}   \put(208,64) {$11$}  \put(238,64) {$13$}
\end{picture}

Consider the complete ladder with some rungs missing. Precisely, let $\ll=\{\l_j\}_{j\ge1}\subset\bn$ be a sequence of positive integers, $\l_1=1$, and
assume that the rungs $\{2n-1, 2n\}_{n\in\ll}$ are present, i.e., the vertices $2n-1$ and $2n$ are incident. Put
$$ \chi_k:=\left\{
             \begin{array}{ll}
               1, & \hbox{$k\in\ll$;} \\
               0, & \hbox{$k\notin\ll$.}
             \end{array}
           \right.  \qquad (\chi_1=1).
$$
We define the canonical basis as in the case of the complete ladder. The subspaces $\cj^\pm=\lc\{h_n^\pm\}_{n\ge1}$ are invariant for $\ag$,
$\ell^2=\cj^+\oplus\cj^-$, and the adjacency operator $\ag$ acts as
\begin{equation*}
\begin{split}
\ag h_1^+ &=h_1^+ +h_2^+, \qquad \ \ \ag h_n^+ =h_{n-1}^+ +\chi_n h_n^+ +h_{n+1}^+, \quad n=2,3,\ldots, \\
\ag h_1^- &=-h_1^- +h_2^-, \qquad \ag h_n^- =h_{n-1}^- -\chi_n h_n^- +h_{n+1}^-, \quad n=2,3,\ldots.
\end{split}
\end{equation*}
So, the canonical form of $\gg$ is
\begin{equation}\label{simlad}
\ag=J^+\bigoplus J^-, \quad J^\pm=J(\{\pm\chi_1, \pm\chi_2,\ldots\}, \{1\}).
\end{equation}

There is no hope computing the spectra of the above matrices $J^\pm$ explicitly. So we will focus here on two particular cases
with periodic structure, and on the sparse case in the next section.

A periodic structure appears when $\ll$ is an arithmetic progression. We consider two simplest examples of such kind.

Let $\ll=\{1,3,5,\ldots\}$ (the infinite linear hexagon $L_\infty$). The terms in \eqref{simlad} are
$$ J^\pm=J(\{\pm1, 0, \pm1, 0, \ldots\}, \{1\}), $$
so we have two $2$-periodic Jacobi matrices.

A detailed computation is carried out for $J^+$. The discriminant \eqref{discr} and the polynomial $\g_2$ in \eqref{weyl} are
$$ \cd^+(\l)=\l(\l-1)-2, \qquad \g_2(\l)=\l(\l-1). $$
So,
$$ \cd^+(\l)=2 \ \sim \ \l_{1,2}^+=\frac{1\pm\sqrt{17}}2\,, \quad \cd^+(\l)=-2 \ \sim \ \l_{3,4}^+=0,1. $$
The essential spectrum is
\begin{equation}\label{hexhon1}
\s_{ess}(J^+)=\Bigl[\frac{1-\sqrt{17}}2\,, 0\Bigr]\bigcup \Bigl[1, \frac{1+\sqrt{17}}2\Bigr].
\end{equation}

The only root of the polynomial $p_2$ is $1$, which is at the edge of the gap. So, there are no eigenvalues, and $\s(J^+)=\s_{ess}(J^+)$.

The computation for $J^-$ is identical. Now the essential spectrum is
\begin{equation}\label{hexhon2}
\s_{ess}(J^-)=\Bigl[\frac{-1-\sqrt{17}}2\,, -1\Bigr]\bigcup \Bigl[0, \frac{-1+\sqrt{17}}2\Bigr],
\end{equation}
and again, there are no eigenvalues of $J^-$. Finally
\begin{equation}\label{hexhon}
\s(\gg)=\s_{ess}(\gg)=\Bigl[\frac{-1-\sqrt{17}}2\,, \frac{1+\sqrt{17}}2\Bigr].
\end{equation}

\smallskip

Let now $\ll=\{1,4,7,10,\ldots\}$ (the infinite linear octagon), so the Jacobi components are
$$ J^\pm=J(\{\pm1, 0, 0, \pm1, 0, 0,\ldots\}, \{1\}), $$
and we have two $3$-periodic Jacobi matrices. Again, a direct computation provides
$$ \cd^+(\l)=\l^3-\l^2-3\l+1, \quad \g_3(\l)=(\l-1)^2(\l+1). $$
The endpoints of the gaps are
\begin{equation*}
\begin{split}
\cd^+ &=2 \ \sim \ \l_1^+=-1,\ \ \l_{2,3}^+=1\pm\sqrt2; \\
\cd^+ &=-2 \ \sim \ \l_4^+=1,\ \ \l_{5,6}^+=\pm\sqrt3,
\end{split}
\end{equation*}
so the essential spectrum is
$$\s_{ess}(J^+)=\bigl[-\sqrt3, -1\bigr]\cup\bigl[1-\sqrt2, 1\bigr]\cup\bigl[\sqrt3, 1+\sqrt2\bigr]. $$
Now the both roots of the polynomial $p_3(\l)=\l^2-\l-1$ lie inside the gaps
$$ \l_7^+=\frac{1+\sqrt5}2\in\bigl(1, \sqrt3\bigr), \quad \l_8^+=\frac{1-\sqrt5}2\in\bigl(-1, 1-\sqrt2\bigr), $$
and $\g(\l_{7,8}^+)>0$. According to the ``rule of signs'', the point $\l_7^+$ is the eigenvalue of $J^+$, and $\l_8^+$ is not.

The computation for $J^-$ is identical. We have
$$\s_{ess}(J^-)=\bigl[-1-\sqrt2,-\sqrt3\bigr]\cup\bigl[-1, -1+\sqrt2\bigr]\cup\bigl[1, \sqrt3\bigr], $$
and the only eigenvalue is $(-1-\sqrt5)/2$. Finally,
\begin{equation}\label{octhon}
\s(\gg)=\bigl[-1-\sqrt2, 1+\sqrt2\bigr]\bigcup \,\Bigl\{\pm\frac{1+\sqrt5}2\Bigr\},
\end{equation}
and the spectrum has multiplicity $2$ on $[1-\sqrt2, \sqrt2-1]$.
\end{example}

\begin{example}\label{squares} ``Squares and cubes''.

Consider the simplest chain of squares

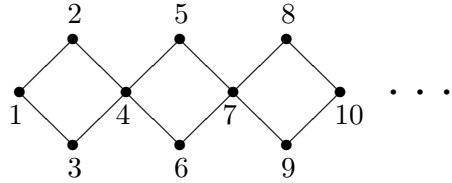
\begin{picture}(250, 100)
\multiput(40,50) (40,0) {4} {\circle* {4}}
\multiput(60,30) (40,0) {3} {\circle* {4}}
\multiput(60,70) (40,0) {3} {\circle* {4}}
\put(40,50) {\line(1,1){20}}  \put(40,50) {\line(1,-1){20}}
\put(60,30) {\line(1,1){20}} \put(80,50) {\line(1,1){20}}
\put(100,30) {\line(1,1){20}} \put(120,50) {\line(1,1){20}} \put(140,30) {\line(1,1){20}}
\put(60,70) {\line(1,-1){20}} \put(80,50) {\line(1,-1){20}} \put(100,70) {\line(1,-1){20}}
\put(120,50) {\line(1,-1){20}} \put(140,70) {\line(1,-1){20}}
\multiput(180,50) (10,0) {3} {\circle* {2}}
\put(36,38) {$1$} \put(76,38) {$4$} \put(116,38) {$7$} \put(158,38) {$10$}
\put(58,18) {$3$} \put(98,18) {$6$} \put(138,18) {$9$}
\put(58,76) {$2$} \put(98,76) {$5$} \put(138,76) {$8$}
\end{picture}

We define the canonical basis as
$$ h_{3k-2}=e_{3k-2}, \quad h_{3k-1}=\frac{e_{3k-1}+e_{3k}}{\sqrt2}\,, \quad h_{3k}=\frac{e_{3k-1}-e_{3k}}{\sqrt2}\,, \quad k=1,2,\ldots. $$
The Hilbert space is decomposed as
$$ \ell^2=\ch_0\oplus\ch_1, \qquad \ch_0:=\lc\{h_{3k}\}_{k\ge1}\,, \quad \ch_1:=\lc\{h_{3k-2}, h_{3k-1}\}_{k\ge1}, $$
both $\ch_0$, $\ch_1$ are $\ag$-invariant, and
\begin{equation*}
\ag=\bo_\infty\bigoplus \sqrt2\,J_0.
\end{equation*}
So the spectrum is $\s(\gg)=\{0_\infty\}\cup [-2\sqrt2, 2\sqrt2]$.

A similar graph below is more delicate.

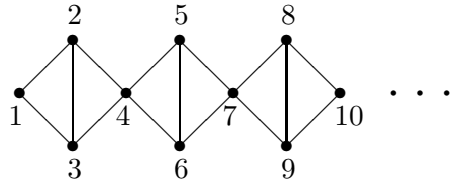
\begin{picture}(250, 100)
\multiput(40,50) (40,0) {4} {\circle* {4}}
\multiput(60,30) (40,0) {3} {\circle* {4}}
\multiput(60,70) (40,0) {3} {\circle* {4}}
\put(40,50) {\line(1,1){20}}  \put(40,50) {\line(1,-1){20}}
\put(60,30) {\line(1,1){20}} \put(80,50) {\line(1,1){20}}
\put(100,30) {\line(1,1){20}} \put(120,50) {\line(1,1){20}} \put(140,30) {\line(1,1){20}}
\put(60,70) {\line(1,-1){20}} \put(80,50) {\line(1,-1){20}} \put(100,70) {\line(1,-1){20}}
\put(120,50) {\line(1,-1){20}} \put(140,70) {\line(1,-1){20}}
\multiput(180,50) (10,0) {3} {\circle* {2}}
\put(36,38) {$1$} \put(76,38) {$4$} \put(116,38) {$7$} \put(158,38) {$10$}
\put(58,18) {$3$} \put(98,18) {$6$} \put(138,18) {$9$}
\put(58,76) {$2$} \put(98,76) {$5$} \put(138,76) {$8$}
\multiput(60,70) (0,-40) {1} {\line(0,-1) {38}}
\multiput(100,70) (0,-40) {1} {\line(0,-1) {38}}
\multiput(140,70) (0,-40) {1} {\line(0,-1) {38}}
\end{picture}

With the same canonical basis and decomposition of $\ell^2$ we have now
\begin{equation}\label{squares}
\ag=-\bi_\infty\bigoplus J(\{0,1,0,1,\ldots\}, \{\sqrt2\}).
\end{equation}
The second term is a $2$-periodic Jacobi matrix. Its discriminant $\cd$ in \eqref{discr} and the polynomial $\g_2$ in \eqref{weyl} are
\begin{equation*}
\cd(\l) =p_3(\l)-p_1^{(1)}(\l)=\frac{\l(\l-1)}2-2, \quad \g_2(\l)=p_3(\l)+p_1^{(1)}(\l)=\frac{\l(\l-1)}2\,.
\end{equation*}
The essential spectrum is a union of two intervals (with the eigenvalue $-1$ of infinite multiplicity on one of them)
\begin{equation}\label{spsq}
\s_{ess}(\gg)=\Bigl[\frac{1-\sqrt{33}}2, 0\Bigr] \bigcup \Bigl[1, \frac{1+\sqrt{33}}2\Bigr].
\end{equation}
The polynomial $p_2=\l/\sqrt2$ has the only root $\l_0=0$ which lies at the edge of the gap, so there are no eigenvalues.
The whole spectrum now is \eqref{spsq} along with the eigenvalue $-1$ of infinite multiplicity.

We might equally well have considered cubes in $\br^3$ in place of squares, with one common vertex for each pair of adjacent cubes.

Define an orthonormal sequence $\{h_n\}_{n\ge1}$ by the recipe
\begin{equation*}
\begin{split}
h_{3k+1} &=e_{3k+1}, \\
h_{3k+2} &=\frac1{\sqrt3}\,\bigl(e_{7k+2}+e_{7k+3}+e_{7k+4}\bigr), \\
h_{3k+3} &=\frac1{\sqrt3}\,\bigl(e_{7k+5}+e_{7k+6}+e_{7k+7}\bigr), \quad k=0,1,\ldots.
\end{split}
\end{equation*}
It is easy to see that the subspace $\ch:=\lc\{h_k\}_{k\ge0}$ is $\ag$-invariant, and the matrix of $\ag\vert\ch$ is a $3$-periodic Jacobi matrix
$$ \ag\vert\ch \sim J(\{0\}, \{\sqrt3,2,\sqrt3,\ldots\}). $$

To complete the above system to an orthonormal basis, we put
\begin{equation*}
\begin{split}
h_{3k+2}^{(1)} &=\frac1{\sqrt6}\,\bigl(e_{7k+2}+e_{7k+3}-2e_{7k+4}\bigr), \quad h_{3k+2}^{(2)} =\frac1{\sqrt2}\,\bigl(e_{7k+2}-e_{7k+3}\bigr),   \\
h_{3k+3}^{(1)} &=\frac1{\sqrt6}\,\bigl(-e_{7k+5}+2e_{7k+6}-e_{7k+7}\bigr), \quad h_{3k+2}^{(2)} =\frac1{\sqrt2}\,\bigl(e_{7k+5}-e_{7k+7}\bigr),
\end{split}
\end{equation*}
$k=0,1,\ldots$. The system thus obtained is complete, and
$$ \ag h_{3k+2}^{(q)}=h_{3k+3}^{(q)}, \quad \ag h_{3k+3}^{(q)}=h_{3k+2}^{(q)}, \quad q=1,2. $$
So, $\ch_k^{(q)}:=\lc\{h_{3k+2}^{(q)}, h_{3k+3}^{(q)}\}$ are invariant subspaces of dimension $2$, and
$$ \ag\vert\ch_k^{(q)}\sim
\begin{bmatrix}
0 & 1 \\
1 & 0
\end{bmatrix}.
$$
Finally,
\begin{equation}\label{cubes}
\ag\sim J(\{0\}, \{\sqrt3,2,\sqrt3,\ldots\})\bigoplus\Bigl(\bigoplus_{n=1}^\infty\Bigr)
\begin{bmatrix}
0 & 1 \\
1 & 0
\end{bmatrix}\Bigr).
\end{equation}

For the Jacobi component we find
\begin{equation*}
\begin{split}
p_4(\l) &=\frac{\l^3-7\l}6\,, \quad p_2^{(1)}(\l)=\frac{\l}2\,, \\
\cd_4(\l) &=\frac{\l^3-10\l}6\,, \quad \g_4(\l)=\frac{\l^3-4\l}6\,,
\end{split}
\end{equation*}
and
\begin{equation}\label{speess}
\s_{ess}(\gg)=[-1-\sqrt7, -2]\cup [1-\sqrt7, \sqrt7-1]\cup [2,1+\sqrt7].
\end{equation}
The denominator in \eqref{weyl}
$$ p_3(\l)=\frac{\l^2-3}{2\sqrt3}\,, \qquad \l_{1,2}=\pm\sqrt3 $$
lie inside the gaps, and the rule of signs says that there are no eigenvalues for $\gg$.
So, $\s(\gg)=\s_{ess}(\gg)$ \eqref{speess} along with the eigenvalues $\pm1$ of infinite multiplicity.
\end{example}

\begin{example}\label{chaincycl} ``A chain of cycles''.

Let $\{n_J\}_{j\ge1}$ be a sequence of positive integers, $n_j\ge2$. Consider a sequence of cycles $\{\bc_{2n_j}\}_{j\ge1}$ of even orders
$|\bc_{2n_j}|=2n_j$. We connect them in a chain in such a way that two adjacent cycles have one vertex in common.

Put
$$ m_j:=2(n_1+\ldots+n_j)-j, \qquad m_0:=0. $$
The vertices of $k$-th cycle $\bc_{2n_k}$ are enumerated with the numbers
$$ \{m_{k-1}+1,\ldots,m_k+1\}, $$
so $\{m_{j-1}+1\}_{j\ge1}$ are the vertices of valency $4$.

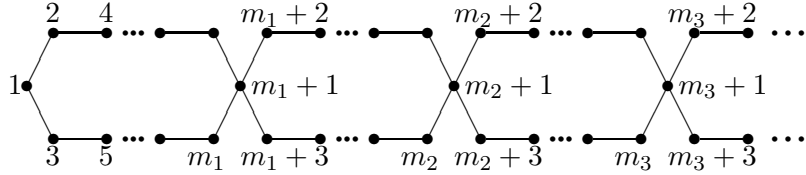
\begin{picture}(350, 120)
\multiput(40,30) (20,0) {14} {\circle* {4}}
\multiput(40,70) (20,0) {14} {\circle* {4}}
\multiput(30,50) (80,0) {4}  {\circle* {4}}
\multiput(67,30) (3,0) {3} {\circle* {2}}
\multiput(67,70) (3,0) {3} {\circle* {2}}
\multiput(42,30) (20,0) {1} {\line(1, 0) {16}}
\multiput(42,70) (20,0) {1} {\line(1, 0) {16}}
\multiput(82,30) (20,0) {1} {\line(1, 0) {16}}
\multiput(82,70) (20,0) {1} {\line(1, 0) {16}}
\multiput(122,30) (20,0) {1} {\line(1, 0) {16}}
\multiput(122,70) (20,0) {1} {\line(1, 0) {16}}
\multiput(162,30) (20,0) {1} {\line(1, 0) {16}}
\multiput(162,70) (20,0) {1} {\line(1, 0) {16}}
\multiput(202,30) (20,0) {1} {\line(1, 0) {16}}
\multiput(202,70) (20,0) {1} {\line(1, 0) {16}}
\multiput(242,30) (20,0) {1} {\line(1, 0) {16}}
\multiput(242,70) (20,0) {1} {\line(1, 0) {16}}
\multiput(282,30) (20,0) {1} {\line(1, 0) {16}}
\multiput(282,70) (20,0) {1} {\line(1, 0) {16}}
\multiput(310,30) (5,0) {3} {\circle* {2}}
\multiput(310,70) (5,0) {3} {\circle* {2}}
\put(30,50) {\line(1,2){10}}  \put(30,50) {\line(1,-2){10}}
\put(110,50) {\line(1,2){10}}  \put(110,50) {\line(1,-2){10}} \put(110,50) {\line(-1,2){10}}  \put(110,50) {\line(-1,-2){10}}
\put(190,50) {\line(1,2){10}}  \put(190,50) {\line(1,-2){10}} \put(190,50) {\line(-1,2){10}}  \put(190,50) {\line(-1,-2){10}}
\put(270,50) {\line(1,2){10}}  \put(270,50) {\line(1,-2){10}} \put(270,50) {\line(-1,2){10}}  \put(270,50) {\line(-1,-2){10}}
\multiput(147,30) (3,0) {3} {\circle* {2}} \multiput(147,70) (3,0) {3} {\circle* {2}} \multiput(227,30) (3,0) {3} {\circle* {2}} \multiput(227,70) (3,0) {3} {\circle* {2}}
\put(23,47) {$1$} \put(114,47) {$m_1+1$} \put(194,47) {$m_2+1$} \put(274,47) {$m_3+1$}
\put(37,74) {$2$} \put(37,20) {$3$} \put(57,74) {$4$} \put(57,20) {$5$}
\put(90,20) {$m_1$} \put(110,74) {$m_1+2$} \put(110,20) {$m_1+3$}
\put(170,20) {$m_2$} \put(190,74) {$m_2+2$} \put(190,20) {$m_2+3$}
\put(250,20) {$m_3$} \put(270,74) {$m_3+2$} \put(270,20) {$m_3+3$}
\end{picture}

The canonical basis $\{h_n\}_{n\ge1}$ is defined by $h_{m_{k}+1}=e_{m_{k}+1}$ for $k=0,1,\ldots$,
\begin{equation*}
\begin{split}
h_{m_k+2j} &=\frac1{\sqrt2}\,\bigl(e_{m_k+2j}+e_{m_k+2j+1}\bigr)\,, \\
h_{m_k+2j+1} &=\frac1{\sqrt2}\,\bigl(e_{m_k+2j}-e_{m_k+2j+1}\bigr)\,, \quad j-1,2,\ldots,n_{k+1}-1.
\end{split}
\end{equation*}
Next, let
\begin{equation*}
\begin{split}
\ch^+ &=\lc\{h_{m_k+1}; h_{m_k+2j}, \ j=1,\ldots,n_{k+1}-1\}_{k=0}^\infty, \\
\ch_k^- &=\lc\{h_{m_k+2j+1}, \ j=1,\ldots,n_{k+1}-1\}, \qquad \ch^- = \bigoplus_{k=1}^\infty \ch_{k-1}^-.
\end{split}
\end{equation*}
We come to decomposition $\ell^2=\ch^+\oplus\ch^-$ on two $\ag$-invariant subspaces.

The action of the adjacency operator $\ag$ can be easily traced. First, put $r_j:=n_1+\ldots+n_j$, $j\in\bn$.
The restriction of $\ag$ on $\ch^+$ is of the form
$$ \ag \vert \ch^+=J(\{0\},\{a_j\}_{j\ge1}), \quad
a_i=\left\{
      \begin{array}{ll}
        \sqrt2, & \hbox{$i=1,r_1,r_1+1,r_2,r_2+1,\ldots$;} \\
        1, & \hbox{othewise.}
      \end{array}
    \right.
$$
Secondly, the restriction $\ag\vert\ch_{k-1}^-=J_{0,n_k-1}$, the discrete Laplacian of order $n_k-1$. So,
\begin{equation}\label{cycles}
\ag\simeq J(\{0\},\{a_n\})\bigoplus\Bigl(\bigoplus_{k=1}^\infty J_{0,n_k-1}\Bigr).
\end{equation}

There is no hope computing the spectrum of the above Jacobi matrix $J$ explicitly. So we will focus here on some particular cases with periodic structure,
and on the sparse case in the next section.

A periodic structure appears when, e.g., all cycles are identical: $n_j=N$, $j=1,2,\ldots$. For $N=2$ we have a simple chain of squares, so assume that $N\ge3$.
The above Jacobi matrix is $N$-periodic with
$$ \{a_1,a_2,\ldots,a_{N-1},a_N\}=\{\sqrt2, \underbrace{1,1\ldots,1}_{N-2}, \sqrt2\}.  $$

The argument in Example \ref{cheb} (see equations \eqref{1stkind1} and \eqref{1stkind2}) provides an explicit expression for the discriminant
$$\cd_N(\l)=T_{N}\Bigl(\frac{\l}2\Bigr)-U_{N-2}\Bigl(\frac{\l}2\Bigr). $$
Since we are unable to solve the equations $\cd_N=\pm2$ for an arbitrary $N$, we restrict ourselves with the cases $N=3,4$.

Let $N=4$ (a chain of octagons). Then $J$ is $4$-periodic Jacobi matrix, and the discriminant and $\g_4$ in \eqref{weyl} are
$$ \cd_4(\l)=\frac{\l^4-6\l^2+4}2\,, \qquad \g_4(\l)=\frac{\l^4-2\l^2}2\,. $$
It is easy to solve $\cd_4=\pm2$ and find the spectral bands and the essential spectrum
$$ \s_{ess}(J)=[-\sqrt6, -2]\cup [-\sqrt2, 0]\cup [0,\sqrt2]\cup [2,\sqrt6]. $$
Note that there is a closed gap at the origin.

The denominator $p_4$ in \eqref{weyl} is now
$$ p_4(\l)=\frac{\l(\l^2-3)}{\sqrt2}\,, $$
with the roots $\l_0=0$, $\l_{1,2}=\pm\sqrt3$. The last two lie inside the gaps, and the rule of signs shows that both of them are the eigenvalues of $J$.
So,
\begin{equation}
\s(J)=[-\sqrt6, -2]\cup [-\sqrt2, \sqrt2]\cup [2,\sqrt6]\cup\{\pm\sqrt3\}.
\end{equation}

There is another part of the spectrum which comes from Laplacians of order $n_k-1=3$ in \eqref{cycles}. By \eqref{spfinlap},
\begin{equation}\label{finlapl}
\s(J_{0,3})=\{0,\pm\sqrt2\},
\end{equation}
each eigenvalue has infinite multiplicity.

The same computation can be carried out for $N=3$ (a chain of hexagons). Now
\begin{equation}
\s(\gg)=\Bigl[\frac{-1-\sqrt{17}}2, \frac{1-\sqrt{17}}2\Bigr]\cup [-1,1]\cup \Bigl[\frac{-1+\sqrt{17}}2, \frac{1+\sqrt{17}}2\Bigr]\cup\{\pm\sqrt2\}
\end{equation}
along with two eigenvalues $\pm1$ of infinite multiplicity.
\end{example}

\subsection{Sparse ladders and chains of cycles}

One can gather some information about the spectrum of the graphs in Examples 4.3 and 4.5 in the situation opposite in a sense to one considered above.
We assume that the graph in question is sparse.

\smallskip

{\bf Example 4.3} (cont.)

Simon's ladder $\gg$ is said to be sparse if the set $\ll$ is sparse in the sense of Example \ref{sparse}
$$ \lim_{i\to\infty}\l_{i+1}-\l_i=+\infty. $$
The adjacency operator $\ag$ \eqref{simlad} is now the orthogonal sum of two sparse Jacobi matrices $J^\pm$ \eqref{simlad}. The sets of all
right limits $RL(J^\pm)$ is available, see \eqref{rlimsim}, so the Last--Simon theorem applies for the description of their essential spectrum.
The nonzero right limits for $J^+$ have the same spectrum. Take $k=0$ and compute the spectrum of $J_{right}^+$ \eqref{rlimsim} by using
the perturbation determinant \eqref{compd} with
$$ T_0=J_0(\bz)=J(\{0\},\{1\})_{n\in\bz}, \quad T=J_{right}^+, \quad T-T_0=(\cdot, e_0)e_0. $$
Note that the resolvent matrix is known \cite{KiSi03}
$$ (T_0-\l)^{-1}=\|r_{ij}(z)\|_{i,j\in\bz}, \quad r_{ij}(z)=\frac{z^{|i-j|}}{z-z^{-1}}\,, \quad \l=z+\frac1z. $$
So,
$$ L\Bigl(z+\frac1z, J_{right}^+\Bigr)=1+r_{00}(z)=\frac{z^2+z-1}{z^2-1}\,, \quad z_{\pm}=\frac{-1\pm\sqrt5}2 $$
are the roots of the perturbation determinant, and one of them, $z_+$, is in $(0,1)$. The spectrum
$$\s\bigl(J_{right}^+\bigr)=[-2,2]\cup\{\sqrt5\}. $$
Similarly,
$$\s\bigl(J_{right}^-\bigr)=[-2,2]\cup\{-\sqrt5\}, $$
and so
\begin{equation}\label{specsparsimlad}
\s_{ess}(\gg)=[-2,2]\cup\{\pm\sqrt5\}.
\end{equation}
We observe here two isolated points $\pm\sqrt5$ of the essential spectrum. As a Jacobi matrix can not have multiple eigenvalues, those two
are accumulation points for the eigenvalues. The endpoints $\pm2$ can also attract some eigenvalues of $\gg$.

The structure of the spectrum on $[-2,2]$ is subtle. Note first, the all the Jacobi matrices in our consideration are finite valued, that is,
the diagonals $\{b_n\}$, $\{a_n\}$ take a finite number of values ($3$ in the latter example). A result of Remling \cite{rem11}, \cite[Theorem 7.4.6]{SiSz}
states that finite valued Jacobi matrices with nontrivial absolutely continuous spectrum are eventually periodic
$$ a_{n+N}=a_n, \quad b_{n+N}=b_n, \quad n\ge n_0. $$
Clearly, finite valued and sparse Jacobi matrices can not be eventually periodic, so their spectra are {\it purely singular}. Hence, the spectra of the sparse
Simon's ladders on $[-2,2]$ are purely singular.

Moreover, assume that a Simon's ladder is strongly sparse:
\begin{equation}\label{strspar}
\limsup_{j\to\infty} \frac{\log(\l_{j+1}-\l_j)}{j}=+\infty.
\end{equation}
Now we can bound the norms of the transfer matrices $\ct_n^\pm$ \eqref{trans}, so the Simon--Stolz theorem applies. Indeed, the main diagonal of $J^\pm$ looks
$$ \{\pm1, \underbrace{0,0\ldots,0}_{\l_2-\l_1-1}, \pm1, \underbrace{0,0\ldots,0}_{\l_3-\l_3-1}, \pm1,\ldots\}. $$
Then
\begin{equation*}
\begin{split}
V_{\pm}(\l) &:=A^{\pm}(\l;a_n,b_n)=
\begin{bmatrix}
\l\mp1 & -1 \\
1 & 0
\end{bmatrix}, \quad n\in\ll, \\
V(\l) &:=A^{\pm}(\l;a_n,b_n)=
\begin{bmatrix}
\l & -1 \\
1 & 0
\end{bmatrix}, \quad n\notin\ll,
\end{split}
\end{equation*}
and so for $n=\l_q$
\begin{equation*}
\begin{split}
\ct_n^{\pm}(\l) &=A^{\pm}(\l;a_n,b_n)A^{\pm}(\l;a_{n-1},b_{n-1})\ldots A^{\pm}(\l;a_1,b_1) \\
&=V_\pm(\l)\,V^{\l_q-\l_{q-1}-1}(\l)\,V_\pm(\l)\,V^{\l_{q-1}-\l_{q-2}-1}(\l)\ldots V_\pm(\l)\,V^{\l_2-\l_1-1}(\l)\,V_\pm(\l).
\end{split}
\end{equation*}
The matrix $V$ is diagonalizable. Indeed, for $-2<\l<2$
$$ |V(\l)-zI_2|=z^2-\l z+1 \ \Rightarrow \ z_\pm(\l)=\frac{\l\pm i\sqrt{4-\l^2}}2\in\bt, $$
the eigenvalues are unimodular. So, if
$$ U(\l)=
\begin{bmatrix}
z_+(\l) & z_-(\l) \\
1 & 1
\end{bmatrix}, \quad  U^{-1}(\l)=\frac1{z_+(\l)-z_-(\l)}\,
\begin{bmatrix}
1 & -z_-(\l) \\
-1 & z_+(\l)
\end{bmatrix},
$$
we have
$$ V^k(\l)=U(\l)\,
\begin{bmatrix}
z_+^k(\l) & 0 \\
0 & z_-^k(\l)
\end{bmatrix}\, U^{-1}(\l), \quad k\in\bn. $$
Hence, for each $k\in\bn$ and $\l\in (-2,2)$ the bounds
$$ \|V^k(\l)\|\le \|U(\l)\|\,\|U^{-1}(\l)\|\le\frac{\|U(\l)\|^2}{\sqrt{4-\l^2}}=C_1(\l) $$
hold. The latter gives the bound for the norm of the transfer matrices
\begin{equation*}
\begin{split}
\|\ct_{\l_{q}}^\pm(\l)\| &\le C_1^{q-1}(\l)\,\|V_\pm(\l)\|^q, \\
\|\ct_j^\pm(\l)\| &\le C_1^{q}(\l)\,\|V_\pm(\l)\|^q, \quad j=\l_q+1,\ldots,\l_{q+1}-1,
\end{split}
\end{equation*}
so
\begin{equation*}
\begin{split}
\sum_{j=\l_q+1}^{\l_{q+1}-1} \|\ct_j^\pm(\l)\|^{-2} &\ge C_2^{-q}(\l)(\l_{q+1}-\l_q-1), \\
\sum_{q=1}^\infty \sum_{j=\l_q+1}^{\l_{q+1}-1} \|\ct_j^\pm(\l)\|^{-2} &\ge \sum_{q=1}^\infty C_2^{-q}(\l)(\l_{q+1}-\l_q-1).
\end{split}
\end{equation*}
The series on the right side diverges in view of \eqref{strspar}. By the Simon--Stolz theorem, the matrices $J^\pm$ have no eigenvalues in $(-2,2)$,
which means that the spectrum of the strongly sparse Simon's ladder is {\it purely singular continuous} on $(-2,2)$.

\smallskip

{\bf Example 4.5} (cont.)

A chain of cycles $\gg$ is said to be sparse if $\{n_k\}_{k\ge1}$ is sparse. The Jacobi component of the adjacency operator $\ag$ \eqref{cycles}
takes the form of the Jacobi matrix $J$ in Example \ref{sparse}. The sets of all right limits $RL(J)$ is available, see \eqref{rlimcycles}, so
the Last--Simon theorem applies for the description of its essential spectrum. The nonzero right limits for $J$ have the same spectrum.
Take $k=0$ and compute the spectrum of $J_{right}$ \eqref{rlimcycles} by using the perturbation determinant \eqref{compd} with $T_0=J_0(\bz)$,
$T=J_{right}$, and
$$ T-T_0=(\cdot, g_1)h_1+(\cdot, g_2)h_2, \quad g_1=h_2=e_1, \quad g_2=h_1=(\sqrt2-1)(e_0+e_2). $$
We have with $\ka=\sqrt2-1$
\begin{equation*}
\begin{split}
L\Bigl(z+\frac1z, J_{right}\Bigr) &=
\begin{vmatrix}
1+\ka\bigl(r_{01}(z)+r_{12}(z)\bigr) & \ka^2\bigl(r_{00}(z)+r_{22}(z)+r_{02}(z)+r_{20}(z)\bigr) \\
r_{11}(z) & 1+\ka\bigl(r_{01}(z)+r_{21}(z)\bigr)
\end{vmatrix} \\
&= \Bigl(1+2\ka\,\frac{z^2}{z^2-1}\Bigr)^2-2\ka^2 z^2\,\frac{z^2+1}{(z^2-1)^2}=\frac{3z^2-1}{z^2-1}\,.
\end{split}
\end{equation*}
Now both roots
$$ z_{\pm}=\pm\frac1{\sqrt3}\in (-1,1), $$
so
\begin{equation}\label{specsparcycles}
\s_{ess}(\gg)=[-2,2]\cup\Bigl\{\pm\frac4{\sqrt3}\Bigr\}.
\end{equation}
Again, two isolated points $\pm\frac4{\sqrt3}$ of the essential spectrum are accumulation points for the eigenvalues.
The endpoints $\pm2$ can also attract some eigenvalues of $\gg$.

Following the line of reasoning in the above example, we see that the spectrum of the sparse chain of cycles is purely singular. Moreover,
for strongly sparse chains of cycles
$$ \limsup_{j\to\infty} \frac{\log n_j}{j}=+\infty $$
the spectrum on $(-2,2)$ is a combination of purely point one (from the finite component) lying on a purely singular continuous spectrum of the
Jacobi component.

\subsection{Toeplitz graphs and a comb graph}

\begin{example}\label{toeplitz}  ``The Toeplitz graphs''.

We say that an infinite graph $\gg$ is a Toeplitz graph if for some enumeration of the vertex set with positive integers the adjacency matrix $\ag$
is Toeplitz, i.e.,
$$ \ag=[\a_{i-k}]_{i,k\ge1}, \qquad \a_j=a_{-j}=0,1. $$
The simplest Toeplitz graph is the infinite path $\bp_\infty$: $A(\bp_\infty)=J_0$, the discrete Laplacian.

Here is another, a bit more sophisticated, Toeplitz graph.

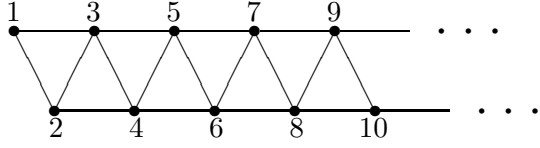
\begin{picture}(250, 90)
\multiput(40,30) (30,0) {5} {\circle* {4}}
\multiput(25,60) (30,0) {5} {\circle* {4}}
\multiput(42,30) (30,0) {5} {\line(1, 0) {26}}
\multiput(27,60) (30,0) {5} {\line(1, 0) {26}}
\put(40,30) {\line(1,2){15}} \put(70,30) {\line(1,2){15}} \put(100,30) {\line(1,2){15}} \put(130,30) {\line(1,2){15}}
\put(25,60) {\line(1,-2){15}} \put(55,60) {\line(1,-2){15}} \put(85,60) {\line(1,-2){15}} \put(115,60) {\line(1,-2){15}} \put(145,60) {\line(1,-2){15}}
\multiput(200,30) (10,0) {3} {\circle* {2}} \multiput(185,60) (10,0) {3} {\circle* {2}}

\put(38,20) {$2$} \put(68,20) {$4$} \put(98,20) {$6$} \put(128,20) {$8$} \put(154,20) {$10$}
\put(22,64) {$1$} \put(52,64) {$3$} \put(82,64) {$5$} \put(112,64) {$7$} \put(142,64) {$9$}
\end{picture}

Clearly, its adjacency matrix is
$$ \ag=[\a_{i-k}]_{i,k\ge1}, \qquad
\a_j=\a_{-j}=\left\{
               \begin{array}{ll}
                 1, & \hbox{$j=1,2$;} \\
                 0, & \hbox{otherwise.}
               \end{array}
             \right.
$$

To find the spectrum of this matrix, we proceed in a standard way. Write the symbol
$$ \p(t)=t^{-2}+t^{-1}+t+t^2=p(\cos\theta), \quad p(x)=4x^2+2x-2. $$
The Toeplitz operator $\ag$ is selfadjoint, so, by \cite[Theorem 1.27]{BoSi99}, the spectrum is
$$ \s(\gg)=p([-1,1])=\Bigl[-\frac94, 4\Bigr]. $$

{\bf Problem 1}. Describe all Toeplitz graphs.
\end{example}

\begin{example}\label{comcomb}  ``A complete comb graph''.

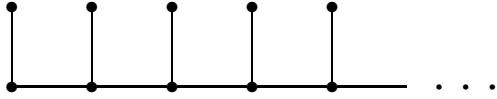
\begin{picture}(200, 90)
\multiput(30,30) (30,0) {5} {\circle* {4}}
\multiput(30,60) (30,0) {5} {\circle* {4}}
\multiput(32,30) (30,0) {5} {\line(1, 0) {26}}
\multiput(190,30) (10,0) {3} {\circle* {2}}
\multiput(30,58) (0,-30) {1} {\line(0,-1) {26}} \multiput(60,58) (0,-30) {1} {\line(0,-1) {26}} \multiput(90,58) (0,-30) {1} {\line(0,-1) {26}}
\multiput(120,58) (0,-30) {1} {\line(0,-1) {26}} \multiput(150,58) (0,-30) {1} {\line(0,-1) {26}}
\end{picture}

The Hilbert space in question is $\ell^2(\cv(\gg))=\ell^2\oplus\ell^2$, and the adjacency operator in the block form looks
$$ \ag=\begin{bmatrix}
J_0 & I \\
I & 0
\end{bmatrix}.
$$
Let us compute the resolvent $R(\ag,z)=(\ag-z)^{-1}$. Clearly,
$$ \begin{bmatrix}
J_0-z & I \\
I & -zI
\end{bmatrix} \begin{bmatrix}
zI & I \\
I & z-J_0
\end{bmatrix}=\begin{bmatrix}
z(J_0-\z(z))I & 0 \\
0 & z(J_0-\z(z))I
\end{bmatrix}, \quad \z(z)=z-\frac1z, $$
and so
$$ R(\ag,z)=\frac1z\,R(J_0,\z(z))\, \begin{bmatrix}
zI & I \\
I & z-J_0
\end{bmatrix}, \quad z\not=0, \quad \z(z)\notin [-2,2]. $$
Since for $z=0$ we have
$$ \ag^{-1}=\begin{bmatrix}
0 & I \\
I & -J_0
\end{bmatrix}, $$
the spectrum is
$$ \s(\gg)=\z^{(-1)}[-2,2]=[-\sqrt2-1, -\sqrt2+1]\cup [\sqrt2-1, \sqrt2+1]. $$
\end{example}

{\bf Problem 2}. As in Example \ref{simlad}, one can consider the complete comb graph with some teeth missing. Analyze the spectral properties of such graphs.

For the description of the spectra of finite comb graphs with tails see \cite{Gol19}.

\section{Spectra of graphs via Schur complement}
\label{s5}

Our second method relies on the block representation \eqref{adjcoup} for the adjacency matrix $\ag$ of the coupling $\gg=\gg_1+\gg_2$ by means of the
bridge of weight $d$. The simplest case is $|\gg_1|=n$, and $d=1$
\begin{equation}\label{adjschur}
A(\gg)=\begin{bmatrix}
A(\gg_1) & E_1 \\
E_1^*& A(\gg_2)
\end{bmatrix}, \qquad
E_1=\begin{bmatrix}
0 & 0 & 0 & \ldots \\
\vdots & \vdots & \vdots &  \\
0 & 0 & 0 & \ldots \\
1 & 0 & 0 & \ldots
\end{bmatrix}.
\end{equation}
Proposition \ref{pro2} applies with $A_{22}=A(\gg_2)$, $A_{11}=A(\gg_1)$, so for $\l\notin\s(A(\gg_2))$ we have $\l\in\s(\ag)$ if and only if the matrix (finite)
$C_{11}(\l)$ \eqref{schucom} is degenerate.

Denote by $P(\cdot,F)$ the characteristic polynomial of a finite graph $F$
$$ P(\l,F)=\det(\l -A(F)). $$
A key ingredient in the argument is the so-called spectral Green's function
$$ G_1(\l,\gg_2)=\bigl(\l-A(\gg_2)^{-1}\bigr)_{11}, \qquad \l\notin\s(A(\gg_2)). $$
Given a graph $F$ and a set of vertices $V\subset \cv(F)$, under $F\backslash V$ we mean the graph induced by the vertices $\cv(F)\backslash V$.

We come to the following result, see \cite[Theorem 1]{Gol17}.

\begin{proposition}\label{mainschur}
Let $\gg=\gg_1+\gg_2$ be the coupling of a finite graph $\gg_1$ and an arbitrary graph $\gg_2$ by means of the bridge $\{n,n+1\}$, and
let $\l\notin\s(\gg_2)$. The point $\l$ belongs to the spectrum of $\gg$  if and only if it solves the equation
\begin{equation}\label{basiceq}
P(\l,\gg_1)-G_1(\l,\gg_2)P(\l,\gg_1\backslash \{n\})=0.
\end{equation}
\end{proposition}

The result in Proposition \ref{mainschur} is effective as long as both Green's function and characteristic polynomials are available.
In the case of infinite tail attached to a finite graph, we have $\gg_2=\bp_\infty$, and, see \eqref{resfree} (and note the negative sign),
$$ G_1(\l, \bp_\infty)=-r_{11}(z)=z, \quad \l=z+\frac1z. $$
The spectrum $\s(\gg)$ in this case is
$$ \s(\gg)=[-2,2]\cup \s_d(\gg), $$
and the discrete spectrum agrees with the roots of the basic equation \eqref{basiceq}
\begin{equation}\label{basiceq1}
\l\in\s_d(\gg) \Leftrightarrow P(\l,\gg_1)-x P(\l,\gg_1\backslash \{n\})=0, \ \ \l=x+\frac1{x}\,,
\ \ x\in(-1,1).
\end{equation}

Let
$$ \s(\gg_1):=\{\l_1\ge\l_2\ge\ldots\ge\l_n\}, \quad \s(\gg_1\backslash \{n\}):=\{\mu_1\ge\mu_2\ge\ldots\ge\mu_{n-1}\} $$
be the spectra of $\gg_1$ and $\gg_1\backslash \{n\}$, respectively. By the Cauchy interlacing theorem,
\begin{equation}\label{inter}
\l_1\ge\mu_1\ge\l_2\ge\ldots\ge\l_{n-1}\ge\mu_{n-1}\ge\l_n,
\end{equation}
(as a matter of fact, the Perron--Frobenius theorem claims that $\l_1>\mu_1$). A quick analysis of the main equation \eqref{basiceq1},
in view of \eqref{inter}, shows that each multiple eigenvalue of $\gg_1$ off $[-2,2]$ solves \eqref{basiceq1}, and so belongs to $\s_d(\gg)$.
One can rewrite \eqref{basiceq1} for $\l>2$ as
\begin{equation}\label{basiceq1.1}
\frac1{x}=\frac{\l+\sqrt{\l^2-4}}2=\frac{P(\l,\gg_1\backslash \{n\})}{P(\l,\gg_1)}=\bigl(\l-A(\gg_1))^{-1}\bigr)_{nn}=G_n(\l, \gg_1),
\end{equation}
the $n$-th spectral Green's function of $\gg_1$. The function on the right side is monotone decreasing on each interval of regularity.
For instance, if the interlacing in \eqref{inter} is strict, i.e.,
$$ \l_1>\mu_1>\l_2>\ldots>\l_{n-1}>\mu_{n-1}>\l_n, $$
$G_n$ is monotone decreasing on each interval $(-\infty, \l_n)$, $(\l_j, \l_{j-1})$, and $(\l_1,\infty)$. Assume that
$$ \l_1>\mu_1>\l_2>\mu_2>\ldots>\mu_{k-1}>\l_{k}>2>\mu_{k}. $$
Then there is exactly one root of \eqref{basiceq1.1} on each interval $(\l_1,\infty)$, $(\l_{j+1},\l_{j})$, $j=1,\ldots,k-1$, and there is no such root on
$(2,\l_k)$. Next, if
$$ \l_1>\mu_1>\l_2>\mu_2>\ldots>\mu_{k-1}>2>\l_{k}>\mu_{k}, $$
then there is exactly one root of \eqref{basiceq1.1} on each interval $(\l_1,\infty)$, $(\l_{j+1},\l_j)$, $j=1,\ldots,k-2$. The existence of the root
on $(\l_{k}, \l_{k-1})$ depends on whether $G_n(2,\gg_1)>1$ (there is a root) or $G_n(2,\gg_1)\le1$ (there are no roots).

The situation for $\l<-2$ can be analyzed in exactly the same way.

\smallskip

To compute the characteristic polynomials in \eqref{basiceq1}, the following result of Schwenk \cite{Schw74}, \cite[Problem 2.7.9]{CDS80}, proves helpful.

{\bf Theorem} (Schwenk). For a given finite graph $F$ and $v\in\cv(F)$, let $\cc(v)$ denote
the set of all simple cycles $Z$ containing $v$. Then
\begin{equation*}
P(\l,F)=\l P(\l,F\backslash v)-\sum_{v'\sim v}P(\l,F\backslash\{v',v\})-
2\sum_{Z\in\cc(v)}P(\l,F\backslash Z).
\end{equation*}

\begin{example}\label{flower} ``A flower with $n$ petals'' \cite{Gol17}.

In this example $\gg_1$ is composed of $n\ge2$ cycles $\{\bc_j\}_{j=1}^n$, glued together at one
common vertex (root) $\co$. Put $\gg=\gg_1+\bp_\infty$ with the infinite path attached to the root $\co$.
Assume that the cycle $\bc_j$ contains $k_j+1\ge3$ vertices.

For the standard Chebyshev polynomial of the second kind $U_k$ \eqref{cheb} we denote
$$ Q(\l,k):=U_k\Bigl(\frac{\l}2\Bigr)=P(\l,\bp_k), \quad Q(\l):=\prod_{j=1}^n Q(\l,k_j)=P(\l,\gg_1\backslash\co), $$
the characteristic polynomials of the finite path $\bp_k$ and the graph $\gg_1\backslash\co$, respectively.

The result of Schwenk applied to the flower graph gives
\begin{equation*}
P(\l,\gg_1)=\l P(\l,\gg_1\backslash \co)-\sum_{v'\sim \co}P(\l,\gg_1\backslash\{v',\co\})-
2\sum_{Z\in\cc(\co)}P(\l,\gg_1\backslash Z).
\end{equation*}
Now
\begin{equation*}
\begin{split}
P(\l,\gg_1\backslash \co) &=\prod_{j=1}^n P(\l,\bp_{k_j})=\prod_{j=1}^n Q(\l,k_j)=Q(\l), \\
P(\l,\gg_1\backslash\{v',\co\}) &=Q(\l,k_m-1)\,\prod_{j\not= m} Q(\l,k_j)=Q(\l)\,\frac{Q(\l,k_m-1)}{Q(\l,k_m)}\,,
\end{split}
\end{equation*}
and so
$$ \sum_{v'\sim \co}P(\l,\gg_1\backslash\{v',\co\})=2Q(\l)\,\sum_{j=1}^n \frac{Q(\l,k_j-1)}{Q(\l,k_j)} $$
(the factor $2$ arises since each cycle enters twice). Next,
$$ P(\l,\gg_1\backslash Z_m)=\prod_{j\not= m} Q(\l,k_j)=\frac{Q(\l)}{Q(\l,k_m)}\,, $$
so
$$ \sum_{m=1}^n P(\l,\gg_1\backslash Z_m)=Q(\l)\,\sum_{j=1}^n \frac1{Q(\l,k_j)}, $$
and finally,
\begin{equation}\label{charflo}
P(\l,\gg_1)=Q(\l)\Bigl\{\l-2\sum_{j=1}^n\frac{Q(\l,k_j-1)+1}{Q(\l,k_j)}\Bigr\}\,.
\end{equation}

Since $Q\not=0$ off $[-2,2]$, the main equation \eqref{basiceq1} looks
$$ \l-2\sum_{j=1}^n\frac{Q(\l,k_j-1)+1}{Q(\l,k_j)}=x, \quad \l=x+\frac1{x} $$
or
\begin{equation}\label{basiceq2}
2\sum_{j=1}^n\frac{Q(\l,k_j-1)+1}{Q(\l,k_j)}=\frac1{x}\,, \quad x\in(-1,1).
\end{equation}

Let first $x=e^{-t}$, $t>0$. We have
$$ Q(\l,k)=Q\Bigl(x+\frac1{x},k\Bigr)=U_k(\cosh)=\frac{\sinh(k+1)t}{\sinh kt}, $$
so
\footnote{There is a misprint in \cite[formula (1.9)]{Gol17}}

\begin{equation}\label{basiceq3}
\p_1(t):=2\sum_{j=1}^n\frac{\sinh k_jt + \sinh t}{\sinh(k_j+1)t}=e^t.
\end{equation}

Note that the functions
$$ f_s(t):=\frac{\sinh at}{\sinh bt}\,, \quad f_c(t):=\frac{\cosh at}{\cosh bt}\,, \quad 0<a<b $$
are monotone decreasing for $t>0$. The latter can be seen, e.g., from the product expansions
$$ \sinh t=t\,\prod_{k=1}^\infty \Bigl(1+\frac{t^2}{k^2\pi^2}\Bigr), \quad \cosh t=\prod_{k=0}^\infty \Bigl(1+\frac{4t^2}{(2k+1)^2\pi^2}\Bigr), $$
(or by elementary calculus). Hence $\p_1$ in \eqref{basiceq3} is monotone decreasing, vanishing at infinity, and $\p(+0)=2n>1$,
so \eqref{basiceq3} has a unique solution $t_+>0$ with
\begin{equation}\label{dspflo1}
\l_+=2\cosh t_+\in\s_d(\gg).
\end{equation}

Similarly, for $x=-e^{-t}$, $t>0$, the main equation takes the form
\begin{equation}\label{basiceq4}
\p_2(t):=2\sum_{j=1}^n\frac{\sinh k_jt +(-1)^{k_j+1}\sinh t}{\sinh(k_j+1)t}=e^t.
\end{equation}
As
$$ \frac{\sinh kt \pm\sinh t}{\sinh(k+1)t}=\frac{\sinh\frac{k\pm 1}2 t\,\cosh\frac{k\mp 1}2 t}{\sinh\frac{k+1}2 t\,\cosh\frac{k+1}2 t}\,, $$
the function $\p_2$ in \eqref{basiceq4} is monotone decreasing (for whatever parity of $k_j$), vanishing at infinity,
$$ \p_2(+0)\ge 2\sum_{j=1}^n \frac{k_j-1}{k_j+1}>1, $$
and \eqref{basiceq4} has a unique solution $t_->0$, so
\begin{equation}\label{dspflo1}
\l_-=-2\cosh t_-\in\s_d(\gg).
\end{equation}
Thereby, the discrete spectrum consists of two points
\begin{equation}\label{disspec2}
\s_d(\gg)=\{\pm 2\cosh t_{\pm}\}.
\end{equation}
\end{example}

In particular case $n=2$, $k_1=k_2$ (the propeller with equal blades), the discrete spectrum is computed in \cite[Example 3.4]{Gol16}.

Note that the method of Section \ref{s3} does not seem to work properly in this example.

\end{document}